\documentclass{amsart}

\usepackage{amssymb}
\usepackage{eepic}
\usepackage{color}
\usepackage{mathtools,xparse}
\usepackage[linktocpage=true]{hyperref}
\usepackage[utf8x]{inputenc}

\allowdisplaybreaks

\numberwithin{equation}{section}

\newtheorem{theorem}{Theorem}[section]
\newtheorem{lemma}[theorem]{Lemma}

\newtheorem{corollary}[theorem]{Corollary}
\newtheorem{remark}[theorem]{Remark}

\newtheorem{ConjA}{Conjecture/Problem (M. de la Salle)}

\newtheorem{TheoA}{Theorem A}

\newtheorem{TheoB}{Theorem $\mathbf{A'}$}

\newcommand{\Z}{\mathbf{Z}}
\newcommand{\R}{\mathbf{R}}
\newcommand{\T}{\mathbf{T}}
\newcommand{\C}{\mathbf{C}}

\newcommand{\B}{\mathcal{B}}
\newcommand{\SL}{S \hskip-1pt L_n(\R)}

\def\G{\mathrm{G}}
\def\1{\mathbf{1}}

\def\Q{\mathcal{Q}}
\def\M{\mathcal{M}}
\def\A{\mathcal{A}}
\def\RR{\mathcal{R}}

\def\S{\mathcal{S}}
\def\V{\mathrm{\mathcal{L}(G)}}

\newcommand{\dem}{\noindent {\bf Proof. }}
\newcommand{\ske}{\noindent {\bf Sketch of the proof. }}

\newcommand{\demInt}{\noindent {\bf Proof of Theorem \ref{Thm-Interpolation}. }}

\newcommand{\fin}{\hspace*{\fill} $\square$ \vskip0.2cm}

\def\mean{- \hskip-11.2pt \int}

\def\meann{- \hskip-10.7pt \int}

\begin{document}

\null

\vskip-50pt

\null

\begin{center}
{\huge Schur multipliers in \\ Schatten-von Neumann classes}

\vskip15pt

{\sc {Jos\'e M. Conde-Alonso, Adri\'an M. Gonz\'alez-P\'erez \\ Javier Parcet and Eduardo Tablate}}
\end{center}

\title[Schur multipliers]{}
%\title[{\sc Schur multipliers in Schatten-von Neumann classes}]{Schur multipliers \\ in Schatten-von Neumann classes}

%\author{Jos\'e M. Conde-Alonso, Adri\'an M. Gonz\'alez-P\'erez \\ Javier Parcet and Eduardo Tablate}

\maketitle

\null

\vskip-45pt

\null

\begin{center}
{\large {\bf Abstract}}
\end{center}

\vskip-25pt

\null

\begin{abstract}
We establish a rather unexpected and simple criterion for the boundedness of Schur multipliers $S_M$ on Schatten $p$-classes which solves a conjecture proposed by Mikael de la Salle. Given $1 < p < \infty$, a simple form our main result reads for $\R^n \! \times \R^n$ matrices as follows $$\big\| S_{\hskip-1pt M} \hskip-2pt: S_p \to S_p \big\|_{\mathrm{cb}} \lesssim \frac{p^2}{p-1} \hskip-2pt \sum_{|\gamma| \le [\frac{n}{2}] +1} \hskip-2pt \Big\| |x-y|^{|\gamma|} \Big\{ \big| \partial_x^\gamma M(x,y) \big| + \big| \partial_y^\gamma M(x,y) \big| \Big\} \Big\|_\infty.$$ In this form, it is a full matrix (nonToeplitz/nontrigonometric) amplification of the H\"ormander-Mikhlin multiplier theorem, which admits lower fractional differentiability orders $\sigma > \frac{n}{2}$ as well. It trivially includes Arazy's conjecture for $S_p$-multipliers and extends it to $\alpha$-divided differences. It also leads to new Littlewood-Paley characterizations of $S_p$-norms and strong applications in harmonic analysis for nilpotent and high rank simple Lie group algebras. 
\end{abstract}

\addtolength{\parskip}{+1ex}

\vskip20pt

\section*{\bf \large Introduction}

Schur multipliers form a class of linear maps on matrix algebras with profound connections in functional analysis, operator algebras,  geometric group theory and harmonic analysis. Their definition is extremely simple, we consider maps $S_M$ determined by $$S_M(A) = \big( M(j,k) A_{jk} \big)_{jk}$$ for any $n \times n$ matrix $A$ and certain function $M: \{1,2,\ldots,n\}^2 \to \C$. More general index sets correspond to operators $A$ acting on $L_2(\Omega,\mu)$ for some $\sigma$-finite measure space $(\Omega,\mu)$. Schur multipliers have played a key role in landmark results since the mid 20th century. The celebrated Grothendieck's inequality is intimately connected to a striking characterization of the operator boundedness of Schur multipliers \cite{Gr,PisBAMS}. The impact of Schur multipliers in geometric group theory and operator algebras was early recognized by Haagerup. His pioneering work on free groups \cite{H} and the research thereafter on semisimple lattices \cite{DCH,CH} encoded deep geometric properties of these groups in terms of approximation properties for Herz-Schur multipliers on their matrix algebras. More recently, strong rigidity properties of high rank lattices were discovered in the remarkable work of Lafforgue and de la Salle \cite{LdlS} studying $L_p$-approximations. The $L_p$-theory has gained a considerable momentum since then \cite{GJP,JMP1,JMP2,dLdlS,MR,MRX,PRS,PRo}. In a somewhat different direction, Schur multipliers have been the fundamental tool in the investigation of operator-Lipschitz functions. The solution of a longstanding problem going back to Krein was achieved by Potapov and Sukochev in \cite{PS} by unraveling Arazy's stronger conjecture on Schur multipliers of divided differences in Schatten $p$-classes. Our main result gives a strikingly simple criterion, which goes beyond best-known to date $S_p$-estimates for divided differences and Herz-Schur multipliers. It also leads to challenging $L_p$-bounds for Fourier multipliers on Lie group algebras. 

In this paper, we shall focus on Schur multipliers with index set $(\Omega,\mu)$ being the $n$-dimensional Euclidean space equipped with its Lebesgue measure. In the context of our main result, this is not a serious restriction and other groups and measure spaces are investigated in \cite{CGPT}. There exists a very fruitful relation between Fourier multipliers and Schur multipliers with symbols $M(x,y) = m(x-y)$, so called Toeplitz or Herz-Schur. Indeed, let $S_p(\R^n)$ denote the Schatten $p$-class acting on $L_2(\R^n)$ and consider the Fourier and Herz-Schur multipliers associated to $m: \R^n \to \C$ $$\widehat{T_m f}(\xi) = m(\xi) \widehat{f}(\xi),$$ \vskip-15pt $$S_m(A) = \big( m(x-y) A_{xy}),$$ with $S_m$ densely defined on $S_2(\R^n) \cap S_p(\R^n)$. Finding accurate conditions on the symbol $m$ which ensure complete $L_p$-boundedness of these multipliers is a rather difficult problem. Building on previous work of Bo\.zejko and Fendler, it is known from \cite{CS,NR} that both problems coincide for this (amenable) group. More precisely
\begin{equation} \tag{FS} \label{Eq-FS}
\hskip4pt \big\| S_m: S_p(\R^n) \to S_p(\R^n) \big\|_{\mathrm{cb}} \, = \, \big\| T_m: L_p(\R^n) \to L_p(\R^n) \big\|_{\mathrm{cb}}. 
\end{equation}
Here, $\| \,\|_{\mathrm{cb}}$ denotes the completely bounded norm of these maps after endowing the corresponding $L_p$ spaces  with their natural operator space structure \cite{P2}. Complete norms are essential in this isometry. The H\"ormander-Mikhlin fundamental theorem \cite{Ho,Mi} gives a criterion for (complete) $L_p$-boundedness of Fourier multipliers. In combination with \eqref{Eq-FS} we may write it as follows 
\begin{equation} \tag{HMS} \label{Eq-HMS}
\hskip-11pt \big\| S_m \hskip-2pt: S_p(\R^n) \to S_p(\R^n) \big\|_{\mathrm{cb}} \, \lesssim \, \frac{p^2}{p-1} \hskip-2pt \sum_{|\gamma| \le [\frac{n}{2}] +1} \Big\| |\xi|^{|\gamma|} \big| \partial_\xi^\gamma m(\xi) \big| \Big\|_\infty.
\end{equation}
The singularity of $\partial_\xi^\gamma m$ around $0$ for $\gamma \neq 0$ is linked to deep concepts in harmonic analysis with profound applications in theoretical physics, differential geometry and partial differential equations. A Sobolev form admits fractional differentiability orders up to $\frac{n}{2} + \varepsilon$ for any $\varepsilon > 0$. This is optimal since we may not consider less derivatives or larger upper bounds for them. Moreover, \eqref{Eq-HMS} is necessary up to order $\frac{n-1}{2}$ for radial $L_p$-multipliers with arbitrarily large $p < \infty$. A characterization of Fourier $L_p$-multipliers for $p \neq 1,2,\infty$ is simply out of reach. 

NonToeplitz multipliers $M(x,y) \neq m(x-y)$ are no longer related to Fourier multipliers. This class is much less understood and their mapping properties are certainly mysterious. The Grothendieck-Haagerup characterization is the main known result. It claims that $S_M$ is bounded on $\mathcal{B}(L_2(\Omega))$ iff it is completely bounded iff there exists a separable Hilbert space $\mathcal{K}$ and two essentially bounded measurable $\mathcal{K}$-valued functions $x \mapsto u_x$ and $y \mapsto w_y$ satisfying 
$$M(x,y) = \langle u_x, w_y \rangle_\mathcal{K} \quad \mbox{for almost every} \quad x,y \in \Omega.$$
The proof for general $\sigma$-finite spaces $(\Omega,\mu)$ can be found in \cite{LdlS}. Pisier characterized \lq\lq unconditional symbols" in connection with lacunary sets in nonamenable groups \cite{PisAJM}. Pisier and Shlyakhtenko applied their Grothendieck's theorem for operator spaces to characterize Schur multipliers which are completely bounded from the space of compact operators to the trace class \cite{PisS}. Bennett also established in \cite{Bennett} connections between completely summing maps and the mapping properties of Schur multipliers acting on $\mathcal{B}(\ell_p, \ell_q)$. Last but not least, Aleksandrov and Peller characterized Hankel-Schur multipliers in $S_p$ $(0 < p < 1)$ with lacunary power series (among many other results) in a very complete paper \cite{AP}. 

Despite a rich literature on Schur multipliers, it becomes rather limited regarding sufficient conditions for $S_p$-boundedness when $1 < p < \infty$. The main results so far are Marcinkiewicz type conditions \cite{BG,CPSW}, a careful analysis of
unconditionality in Schatten $p$-classes in analogy with $\Lambda(p)$-set theory \cite{Ha}, and the alluded solution of Arazy's conjecture \cite{Ar,PS}. In 2019, during the \'Ecole d'automne \lq\lq Fourier Multipliers on Group Algebras\rq\rq${}$ at Besançon, Mikael de la Salle formulated a vague conjecture in search of a general class of \emph{singular Schur multipliers}. In 2020 he referred to it at UCLA as a dream statement for a purely noncommutative HM theorem \cite{MdlS}:

\vskip5pt

\begin{ConjA}
In the spirit of inequality \eqref{Eq-HMS}, find a regularity condition on $($\hskip-1pt nonToeplitz$)$ $C^{[n/2]+1}$-symbols $M: \R^n \! \times \R^n \to \C$ outside the diagonal with a controlled explosion on the diagonal which implies the $($\hskip-1pt complete$)$ $S_p$-boundedness of the associated Schur multiplier $S_M$ for $1 < p < \infty$. 
\end{ConjA}

\vskip-3pt 

Our approach towards the solution is based on a rather elementary observation which led us to think on a possible relation between the theory of Schur multipliers and that of pseudodifferential operators. Let us momentarily argue with the  group $\Z$ to avoid topological considerations. We know from Fourier-Schur transference \eqref{Eq-FS} (valid for any amenable group) that every Toeplitz symbol $M(j,k) = m(j-k)$ gives rise to a Schur multiplier which is identified with the Fourier multiplier $T_m$ on the torus $\mathbf{T}$. On the other hand, if $\bullet$ and $\cdot$ respectively stand for the Schur and ordinary matrix products on $\mathcal{B}(\ell_2(\Z))$, elementary nonToeplitz symbols of the form $M_\alpha(j,k) = \alpha(j)$ and $M_\beta(j,k) = \beta(k)$ satisfy the simple identities $$M_\alpha \bullet A = \mathrm{diag}(\alpha) \cdot A \quad \mathrm{and} \quad M_\beta \bullet A = A \cdot \mathrm{diag}(\beta)$$ for the diagonal matrices with entries $(\alpha(j))_{j \in \Z}$ and $(\beta(k))_{k \in \Z}$. In other words, we have identified \lq\lq Fourier multipliers\rq\rq${}$ and (left/right) \lq\lq pointwise multipliers\rq\rq${}$ within the class of Schur multipliers. Then, recalling that every symbol $M$ can be rewritten as follows  
\begin{eqnarray*}
M(j,k) \!\! & = & \!\! M_r(k-j,k) \ = \ \sum_{\ell \ge 1} m_{r\ell}(j-k) \beta_\ell(k), \\
M(j,k) \!\! & = & \!\! M_c(j,j-k) \ = \ \sum_{\ell \ge 1} \alpha_\ell(j) m_{c\ell}(j-k), 
\end{eqnarray*}
Schur multipliers are combinations of Fourier and left/right pointwise multipliers. Of course, this arises from the well-known isomorphism $\mathcal{B}(\ell_2(\Z)) \simeq \ell_\infty(\Z) \rtimes \Z$ and a transpose form of it. The formal analog of the above reasoning in Euclidean spaces leads to similar conclusions. It is thus conceivable to think of Schur multipliers on Euclidean spaces as matrix forms of pseudodifferential operators \cite{Ho2,Ta}. In particular, our viewpoint might look to be rather close to the theory developed in \cite{GJP2} and one should expect a regularity condition in terms of infinitely many mixed $\partial_x \partial_y$ derivatives of the symbol $M$. Surprisingly, our main result only involves a finite number of $\partial_x$ and $\partial_y$ derivatives (no mixtures) and allows a singular behavior at the diagonal $x=y$, as predicted by de la Salle. 

\begin{TheoA} \label{ThmA1}
If $M \in C^{[n/2]+1}(\R^{2n} \setminus \{x=y\})$ and $1 < p < \infty$
$$\big\| S_M \hskip-2pt: S_p(\R^n) \to S_p(\R^n) \big\| \, \lesssim \, \frac{p^2}{p-1} \big|\big|\big| M \big| \big| \big|_{\mathrm{HMS}},$$ where we set $||| M |||_\mathrm{HMS} := \displaystyle \sum_{|\gamma| \le [\frac{n}{2}] +1} \Big\| |x-y|^{|\gamma|} \Big\{ \big| \partial_x^\gamma M(x,y) \big| + \big| \partial_y^\gamma M(x,y) \big| \Big\} \Big\|_\infty.$
\end{TheoA} 

\null

\vskip-38pt

\null

\includegraphics[scale=0.6]{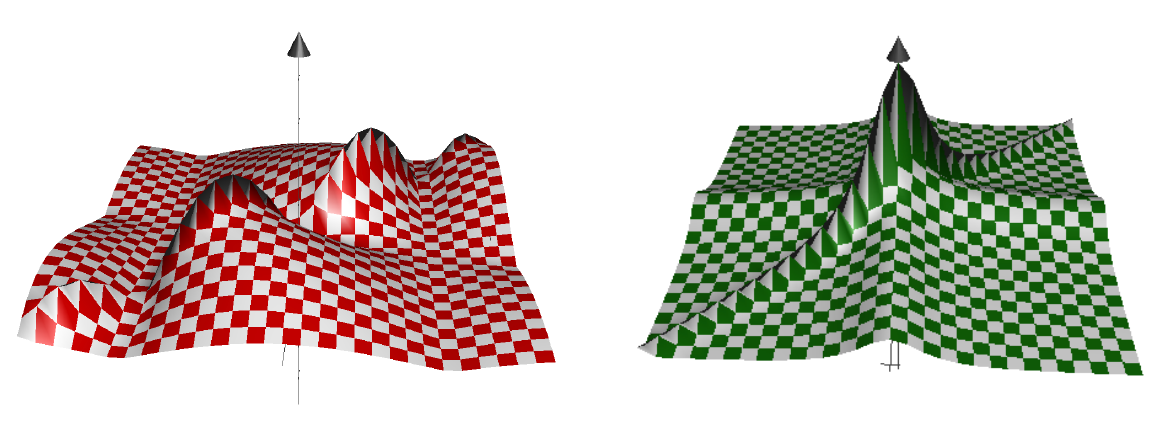}

\vskip-30pt

\null

\begin{center}
{\scriptsize Examples of NonToeplitz H\"ormander-Mikhlin-Schur multipliers in $\R \times \R$} \\ [-3pt]
{\scriptsize Any Toeplitz symbol would be forced to be constant at $x=y+\alpha$ for all $\alpha \in \R$, unlike above}
\end{center}

%\vskip-3pt

Theorem A introduces a new class of $S_p$-bounded multipliers which we shall refer to as H\"ormander-Mikhlin-Schur multipliers. Observe that our condition for Toeplitz symbols $M(x,y) = m(x-y)$ reduces to the classical HM condition on $m$. According to \eqref{Eq-FS}+\eqref{Eq-HMS}, Theorem A is the full matrix (nonToeplitz/nontrigonometric) form of the H\"ormander-Mikhlin multiplier theorem. It admits even weaker regularity assumptions (up to fractional differentiability orders $\sigma > n/2$) expressed in terms of Sobolev norms and even weaker assumptions \`a la Calder\'on-Torchinsky for specific values of $p$, more details in Theorems $\mathrm{A'}$ and Corollary \ref{Cor-CT} below. Moreover, the statement above automatically implies cb-boundedness by \cite[Theorem 1.18]{LdlS}.

The proof is more ingenious than technical. The key novelty is to control the Schur multiplier by two (not just one) \lq\lq twisted Fourier" multipliers. Then, we refine noncommutative Calder\'on-Zygmund methods to prove two different endpoint inequalities for the twisted multipliers. Using a canonical form of BMO for matrix algebras |introduced in this paper and for which we reprove after \cite{JM} the expected complex interpolation behavior with the $S_p$-scale| both inequalities together yield a BMO endpoint inequality which implies the statement. Note that the twisted multipliers satisfy one endpoint inequality (row or column form) but not the other! With hindsight, this illustrates that singular Schur multipliers are slightly better behaved than operator-valued Calder\'on-Zygmund operators. 

The class of HMS multipliers |enlarged under Sobolev regularity| leads to strong applications in operator algebras and harmonic analysis. Arazy's conjecture for divided differences \cite{PS} now follows as a particular example. Moreover, we prove that $\alpha$-divided differences of $\alpha$-H\"older functions are also $S_p$-bounded as long as we have $|1/p - 1/2| < \min\{\alpha,\frac12\}$. It shows a phase transition at $\alpha = 1/2$. This is a very satisfactory result which goes beyond Arazy's conjecture. Second, we obtain strikingly simple $L_p$-boundedness criteria for Fourier multipliers on Lie group von Neumann algebras, in terms of Lie derivatives of the Fourier symbol. This includes nilpotent Lie groups with subRiemannian metrics, high rank simple Lie groups with exponentially growing metrics and a local HM theorem for arbitrary Lie groups with the Riemannian metric. Altogether this goes far beyond \cite{PRS} and we get optimal regularity conditions. On the other hand, the main results in the cocycle-based approach from \cite{JMP1,JMP2} can be deduced and refined for nonorthogonal cocycles as well as for nonunimodular groups. In addition, Theorem A yields new matrix forms of the Littlewood-Paley theorem and admits a discrete formulation. A proof will be given for those applications in the scope of matrix algebras. The rest involve generalized forms of Theorem A and additional techniques postponed to \cite{CGPT}. 

\section{\bf \large Auxiliary results}

Let us begin with a brief discussion on Schur multipliers over a general $\sigma$-finite measure space $(\Omega,\mu)$. We just collect a few basic but necessary notions for what follows. A much more detailed exposition can be found in Section 1 of \cite{LdlS}. Next we introduce a BMO space for matrix algebras and study interpolation with respect to the $S_p$-scale. The interpolation result was already known for an equivalent BMO seminorm constructed in \cite{JM} for more general von Neumann algebras. Our contribution here is limited to introducing a canonical definition of BMO for matrix algebras and to giving a simple self-contained proof of the interpolation theorem.

\subsection{Schur multipliers} \label{Sect-Schur}

Given a $\sigma$-finite measure space $(\Omega,\mu)$ |in this paper we shall be mostly concerned with $(\Omega,\mu)$ being an $n$-dimensional Euclidean space with its Lebesgue measure| let $\B(L_2(\Omega))$ be the matrix algebra of bounded linear operators  on $L_2(\Omega)$ and let $S_\infty(\Omega)$ be the subalgebra of compact operators equipped with the operator norm. If $p \ge 1$, the Schatten $p$-class $S_p(\Omega)$ is the subspace of bounded linear operators on $L_2(\Omega)$ with finite $p$-norm $$S_p(\Omega) = \Big\{ A \in \mathcal{B}(L_2(\Omega)): \|A\|_p = \mathrm{tr} (|A|^p)^{\frac{1}{p}} < \infty \Big\}.$$ More generally, we may consider noncommutative $L_p$-spaces over a von Neumann algebra $\M$ equipped with a normal semifinite faithful trace $\tau$ \cite{PX2}. A linear map $\Lambda: L_p(\M,\tau) \to L_p(\M,\tau)$ is called completely bounded when $\Lambda \otimes \mathrm{id}$ extends to a bounded operator on $L_p(\M \bar\otimes \B(\mathcal{K}), \tau \otimes \mathrm{tr})$ for any Hilbert space $\mathcal{K}$. If $M_n$ denotes the algebra of $n \times n$ matrices, the cb-norm of $\Lambda$ is $$\big\| \Lambda: L_p(\M,\tau) \to L_p(\M,\tau) \big\|_{\mathrm{cb}} = \sup_{n \ge 1} \big\| \mathrm{id} \otimes T: L_p(M_n(\M)) \to L_p(M_n(\M)) \big\|.$$ 

Informally, a Schur multiplier associated to a measurable function $M: \Omega \times \Omega \to \C$ is the map $S_M$ sending an operator $A: L_2(\Omega) \to L_2(\Omega)$ with matrix representation $(A_{xy})_{x,y \in \Omega}$ to the operator represented by $(M(x,y) A_{xy})_{x,y \in \Omega}$. However, not every such operator admits a kernel or a matrix representation and we must give a more precise definition. When $p=2$, every $A \in S_2(\Omega)$ admits a matrix representation and $S_2(\Omega)$-bounded Schur multipliers $S_M$ are in one-to-one correspondence with symbols $M \in L_\infty(\Omega \times \Omega)$. In addition, we  shall say that $S_M: S_p(\Omega) \to S_p(\Omega)$ is (completely) bounded for other values of $1 \le p \le \infty$ when its restriction to $S_2(\Omega) \cap S_p(\Omega)$ extends to a (completely) bounded map. By density in Schatten $p$-classes, the extension is necessarily unique. We shall be using two important properties of Schur multipliers:
\begin{itemize}
\item[i)] If $(\Omega,\mu)$ has no atoms and $1 \le p \le \infty$, then $$\hskip20pt \big\| S_M: S_p(\Omega) \to S_p(\Omega) \big\|_{\mathrm{cb}} = \big\| S_M: S_p(\Omega) \to S_p(\Omega) \big\|.$$

\item[ii)] Assume that $\Omega$ is locally compact and $\mu$ is a $\sigma$-finite Radon measure on $\Omega$. Let $1 \le p \le \infty$, then the following holds for any $M: \Omega \times \Omega \to \C$ continuous
$$\hskip30pt \big\| S_M: S_p(\Omega) \to S_p(\Omega) \big\| = \sup_{\Sigma \in \mathcal{J}(\Omega,\mu)} \big\| S_{M_{\mid_{\Sigma \times \Sigma}}}: S_p(\Sigma) \to S_p(\Sigma) \big\|,$$ with $\mathcal{J}(\Omega,\mu)$ the finite sets in the support of $\mu$ and $S_p(\Sigma)$ constructed with the counting measure on $\Sigma$. A similar result holds as well for the cb-norm.    
\end{itemize}
The proofs of both results above appear in \cite[Theorems 1.18 and 1.19]{LdlS} respectively.

\subsection{Matrix BMO}

Noncommutative BMO spaces were introduced by Pisier and Xu in \cite{PX} for noncommutative martingales and further investigated by Junge and Mei in other contexts \cite{JM,Mei07}. In this paragraph we introduce a new form of noncommutative BMO for matrix algebras and reprove a well-known interpolation identity. This requires two auxiliary forms of operator-valued BMO due to Mei.  

\noindent \textbf{Mei's BMO spaces.} Let 
\begin{eqnarray*}
\mathcal{A} \!\!\!\! & = & \!\!\!\! L_\infty(\T^n) \bar\otimes \mathcal{B}(L_2(\R^n)), \\
\mathcal{R} \!\!\!\! & = & \!\!\!\! L_\infty(\R^n) \bar\otimes \mathcal{B}(L_2(\R^n)),
\end{eqnarray*}
be the von Neumann algebra formed by essentially bounded $\B(L_2(\R^n))$-valued weakly measurable functions on $\T^n$ or $\R^n$. When $p < \infty$, the noncommutative $L_p$ space $L_p(\RR)$ is the $S_p$-valued Euclidean $L_p$-space $L_p(\R^n; S_p(\R^n))$ equipped with the Lebesgue measure. Mei investigated in \cite{Mei07} the seminorm 
$$\|f\|_{\mathrm{BMO}_\RR^c} = \sup_{Q \in \Q} \Big\| \Big( \mean_Q \big| f(z) - f_Q \big|^2 dz \Big)^\frac12 \Big\|_{\mathcal{B}(L_2(\R^n))}.$$ Here $\Q$ is the set of $n$-dimensional balls (equivalently cubes with sides parallel to the axes) and $f_Q$ stands for the $Q$-mean of an operator-valued function $f$. Mei defined $\mathrm{BMO}_\RR^c$ as the set of those elements in the space $\mathcal{J}_c(\nu) = \B(L_2(\R^n)) \bar\otimes L_2^c(\nu)$ having finite $\mathrm{BMO}_\RR^c$-seminorm, where $\nu$ is a probability measure in $\R^n$ irrelevant for our goals. A precise definition of $\mathcal{J}_c(\nu)$ can be found in \cite{Mei07}.  He also proved that $\mathrm{BMO}_\RR^c$ is a dual space, whose predual space is an operator-valued Hardy space. The space $\mathrm{BMO}_\A^c$ is defined similarly with balls/cubes in the torus $\T^n$.

\noindent \textbf{Matrix BMO space.} Consider the unitary $$uf(x,z) = \exp \big( 2\pi i \langle x,z \rangle \big) f(x,z) \quad \mbox{on $L_2(\R^n \times \R^n)$}.$$ The $*$-homomorphism $\pi: \B(L_2(\R^n)) \to L_\infty(\RR)$ by $u$-conjugation is
$$\pi (A) := u \big( \1 \otimes A \big) u^* = \Big( \exp \big( 2\pi i \langle \, \cdot \, , x-y \rangle \big) A_{xy} \Big)$$ where the last expression holds for the weak-$*$ dense subspace $S_2(\R^n) \subset \B(L_2(\R^n))$ and allows us to interpret $\pi(A)$ as an operator-valued function. Then, we define the matrix BMO seminorm 
$$\|A\|_{\mathrm{BMO}} = \max \big\{ \|A\|_{\mathrm{BMO}_r}, \|A\|_{\mathrm{BMO}_c} \big\},$$ 
$$\|A\|_{\mathrm{BMO}_c} = \big\| \pi(A) \big\|_{\mathrm{BMO}_\RR^c} \quad \mbox{and} \quad \|A\|_{\mathrm{BMO}_r} = \|A^*\|_{\mathrm{BMO}_c}.$$
Let $\mathrm{BMO}_c$ be the weak-$*$ closure of $\pi(\B(L_2(\R^n)))$ in the dual space $\mathrm{BMO}_\RR^c$ and equip it with the operator space structure which provides $M_m [\mathrm{BMO}_c]$ with the above $\mathrm{BMO}_c$-seminorm in terms of the amplified $*$-homomorphism $\pi_m = \mathrm{id}_{M_m} \otimes \pi$ instead of $\pi$. Note that $\pi(S_2(\R^n))$ sits in $\mathcal{J}_c(\nu)$ since $$\Big\| \int_{\R^n} \big| \pi(A) \big|^2 d\nu \Big\|_\infty^\frac12 \le \Big\| \int_{\R^n} \big| \pi(A) \big|^2 d\nu \Big\|_1^\frac12 = \Big( \int_{\R^n} \mathrm{tr} \big( \pi |A|^2 \big) d\nu \Big)^\frac12 = \|A\|_{S_2(\R^n)}.$$ Here it is crucial that $\nu$ is a probability measure and it turns out that the couple $(\mathrm{BMO}_c, S_2(\R^n))$ is compatible for complex interpolation. A similar definition also applies for $\mathrm{BMO}_r$ and we set as usual $\mathrm{BMO} = \mathrm{BMO}_r \cap \mathrm{BMO}_c$. In the operator space context, both $(\mathrm{BMO}_r, S_2(\R^n))$ and $(\mathrm{BMO}_c, S_2(\R^n))$ are compatible since all operator spaces sit in $\mathcal{J}_r(\nu) + \mathcal{J}_c(\nu) = \B(L_2(\R^n)) \bar\otimes L_2^{r+c}(\nu)$. We refer to \cite{JLMX} for precise definitions of $L_2^r, L_2^c$ and $L_2^{r+c}$. 

\begin{theorem} \label{Thm-Interpolation}
Complex interpolation gives 
$$\big[ \mathrm{BMO}, S_2(\R^n) \big]_{\frac{2}{p}} \simeq_{\mathrm{cb}} \, S_p(\R^n) \quad \mbox{for} \quad 2 \le p < \infty.$$
The equivalence constant in the complete isomorphism above $c_p \approx p$ as $p \to \infty$.
\end{theorem}

The idea of the proof is elementary, we will just approximate the norm of $A$ in $[\mathrm{BMO}, S_2(\R^n)]_{2/p}$ by the norm in $[\mathrm{BMO}_\A, L_2^\circ(\A)]_{2/p}$ of a mean-zero sequence $\pi_k \circ J_k(A)$ for which the interpolation result holds easily. A careful justification of this approximation process is a bit technical. Our choice of the operator-valued probability space $\A$ makes it possible to find the sequence $\pi_k \circ J_k(A)$, satisfying certain conditions which are crucial along the proof. We start with a well-known approximation result from complex interpolation theory. 

\begin{lemma} \label{Lemma-Interp-Density}
Let $\mathcal{S} = \{z \in \C: 0 < \mathrm{Re}(z) < 1\}$ and let $\partial_j = \{z \in \C: \mathrm{Re}(z) = j\}$ for $j=0,1$. Let $\mathcal{F}(\mathrm{X}_0, \mathrm{X}_1)$ be the space of holomorphic functions $\Psi: \mathcal{S} \to \mathrm{X}_0 + \mathrm{X}_1$ which extend to a continuous function in the closure of $\mathcal{S}$ and whose restrictions $\psi_j: \partial_j \to \mathrm{X}_j$ to the boundaries $\partial_j$ belong to $C_0(\partial_j; \mathrm{X}_j)$. Then, if $A \in S_2(\R^n)$ and $2 < p < \infty$ and $\varepsilon > 0$, \hskip-1pt there exists $\Psi_{p,\varepsilon} \in \mathcal{F}(\mathrm{BMO}, S_2(\R^n))$ satisfying$\hskip1pt :$
\begin{itemize}
\item[i)] $\Psi_{p,\varepsilon}(2/p) = A$ and $\Psi_{p,\varepsilon}(\partial_0 \cup \mathcal{S} \cup \partial_1) \subset S_2(\R^n)$. 

\vskip3pt

\item[ii)] $\Psi_{p,\varepsilon}$ is a linear combination of functions of the form $$\hskip20pt \exp(\delta z^2) \sum_{j=1}^m A_j \exp(\lambda_j z) \quad \mbox{with} \quad (A_j, \lambda_j, \delta) \in S_2(\R^n) \times \R \times \R_+.$$ 

\item[iii)] $\displaystyle \sup_{s \in \R} \max \Big\{ \big\| \Psi_{p,\varepsilon}(is) \big\|_{\mathrm{BMO}}, \big\| \Psi_{p,\varepsilon}(1+is) \big\|_{S_2(\R^n)} \Big\} \le (1+\varepsilon) \big\| A \big\|_{[\mathrm{BMO}, S_2(\R^n)]_{\frac2p}}$.
\end{itemize}
\end{lemma}

\dem According to the complex interpolation method, we know there exists $\Phi_{p,\varepsilon}$ in $\mathcal{F}(\mathrm{BMO}, S_2(\R^n))$ satisfying $\Phi_{p,\varepsilon}(2/p)=A$ and property iii). On the other hand by \cite[Lemma 4.2.3]{BL} we know that elements of $\mathcal{F}(\mathrm{BMO}, S_2(\R^n))$ taking values in $\mathrm{BMO} \cap S_2(\R^n) = S_2(\R^n)$ and satisfying property ii) are norm dense, with the norm given by the left hand side of iii). This completes the proof. \fin

Next, we consider periodic perturbations of $\pi(A)$ as follows. Given $k \in \Z_+$ and $x = (x_1, x_2, \ldots, x_n) \in \R^n$, define $m_k(x) = (m_k(x)_1, m_k(x)_2, \ldots, m_k(x)_n) \in \Z^n$ by $8^{-k} m_k(x)_j \le x_j < 8^{-k} (m_k(x)_j+1)$ for $1 \le j \le n$, so that $m_k(x)_j = [8^kx_j]$ and we set $p_k(x) = 8^{-k} m_k(x)$. Given $s \in \T^n$, construct the map $\pi_k: S_2(\R^n) \to L_2(\A)$ as follows $$\pi_k(A)(s) = \Big( e^{2\pi i \langle s, m_k(x) - m_k(y) \rangle} A_{xy} \Big).$$ Identifying $\T^n \simeq [-\frac12, \frac12]^n$, we also define for $\xi \in [-\frac12 8^k, \frac12 8^k]^n$ $$ \widetilde{\pi}_k(A) (\xi) = \Big( e^{2\pi i \langle \xi, p_k(x) - p_k(y) \rangle} A_{xy} \Big) = \pi_k(A)(8^{-k}\xi).$$ Next, we consider the Toeplitz projection 
$$J_k(A) = \Big( \chi_{|x-y| \ge 2^{-k}} A_{xy} \Big) \quad \mbox{for} \quad A \in S_2(\R^n).$$ In other words, $J_k$ is the Herz-Schur multiplier which removes the diagonal strip of width $2^{-k}$ in the Euclidean metric. We could have used a smooth version of it with better boundedness properties, but this will be enough for our purposes. It is clear that $\pi_k \circ J_k = J_k \circ \pi_k$ and $\pi_k \circ J_k(A)$ is mean-zero in $\T^n$ for any $A \in S_2(\R^n)$
$$\int_{\T^n} \pi_k \circ J_k(A)(s) ds = 0.$$ Indeed, we get nonvanishing frequencies $m_k(x)-m_k(y) \in \Z \setminus \{0\}$ for $J_k(A)_{xy} \neq 0$.

\begin{lemma} \label{Lemma-PikPi}
Given $A \in S_2(\R^n)$, we have $$\big\| \pi_k \circ J_k(A) \big\|_{\mathrm{BMO}_\A} \le \big\| \pi \circ J_k(A) \big\|_{\mathrm{BMO}_\RR} + C_n 2^{-k} \big\| J_k(A) \big\|_{S_2(\R^n)}.$$ 
\end{lemma}

\dem Given a cube $S \subset [-\frac12,\frac12]^n$, let $Q = 8^k S$ and observe that 
\begin{eqnarray*}
\lefteqn{\hskip-30pt \meann_S \Big| \pi_k J_k(A) (u) - \meann_S \pi_k J_k(A) (v) dv \Big|^2 du} \\ \hskip30pt \!\!\!\! & = & \!\!\!\! \meann_Q \Big| \widetilde{\pi}_k J_k(A) (\xi) - \meann_Q \widetilde{\pi}_k J_k(A) (\eta) d\eta \Big|^2 d\xi \\ \hskip30pt \!\!\!\! & = & \!\!\!\! \underbrace{\meann_Q \Big| \Big( \big\{ e^{2\pi i \langle \xi, p_k(x) - p_k(y) \rangle} - \meann_Q e^{2 \pi i \langle \eta, p_k(x) - p_k(y) \rangle} d\eta \big\} J_k(A)_{xy} \Big) \Big|^2 d\xi}_{M\hskip-1pt O_Q(\widetilde{\pi}_kJ_k(A))^2}.
\end{eqnarray*}
Then, we may write 
\begin{eqnarray*}
\big\| \pi_k \circ J_k (A) \big\|_{\mathrm{BMO}_\A^c} \!\!\!\! & = & \!\!\!\! \sup_{\begin{subarray}{c} Q \in \Q \\ Q \subset [-\frac12 8^k, \frac12 8^k]^n \end{subarray}} \big\| M\hskip-1pt O_Q(\widetilde{\pi}_kJ_k(A)) \big\|_{\B(L_2(\R^n))} \\ \!\!\!\! & \le & \!\!\!\! \sup_{Q \ \mathrm{small}} \big\| M\hskip-1pt O_Q(\widetilde{\pi}_kJ_k(A)) \big\| + \sup_{Q \ \mathrm{large}} \big\| M\hskip-1pt O_Q(\widetilde{\pi}_kJ_k(A)) \big\|,
\end{eqnarray*}
where we set $Q$ to be small when $\mathrm{diam} (Q) \le 4^k$ and large otherwise. If the above supremum is nearly attained at a small cube $Q$ with center $c_Q$, we let $Q_0 = Q - c_Q$ be the translated cube centered at $0$. Recall that $\widetilde{\pi}_k(B)(\xi) = u_k(\xi) (\1 \otimes B) u_k^*(\xi)$ for certain unitary $u_k(\xi)$. This implies that $$M\hskip-1pt O_Q(\widetilde{\pi}_kJ_k(A)) = M\hskip-1pt O_{Q_0}(\widetilde{\pi}_kJ_k(A)(\hskip1pt \cdot + c_Q)) = u_k(c_Q) M\hskip-1pt O_{Q_0}(\widetilde{\pi}_kJ_k(A)) u_k(c_Q)^*.$$ Therefore, the supremum is also nearly attained at $Q_0$ and 
\begin{eqnarray*}
\lefteqn{\hskip-5pt \big\| M\hskip-1pt O_{Q_0}(\widetilde{\pi}_kJ_k(A)) \big\|} \\ [5pt] \!\!\! & \le & \!\!\! \big\| M\hskip-1pt O_{Q_0}(\pi J_k(A)) \big\| + 2 \sup_{\xi \in Q_0} \big\| (\widetilde{\pi}_k - \pi ) J_k(A)(\xi) \big\|_{S_2(\R^n)} \\ \!\!\! & \le & \!\!\! \big\| \pi \circ J_k(A) \big\|_{\mathrm{BMO}_\RR^c} + 2 \sup_{\begin{subarray}{c} \xi \in Q_0 \\ x,y \in \R^n \end{subarray}} \big| 1- e^{2\pi i \langle \xi, (x-p_k(x)) - (y-p_k(y)) \rangle} \big| \big\| J_k(A)\big\|_{S_2(\R^n)}. 
\end{eqnarray*}
The latter supremum is bounded by $C_n \mathrm{diam}(Q_0) 8^{-k} \|J_k(A)\|_2 \le C_n 2^{-k} \|J_k(A)\|_2$ which gives the desired estimate. Assume now that the supremum defining the $\mathrm{BMO}_\A^c$ seminorm is nearly attained at a large cube $R$. Then we observe $$\big\| M \hskip-1pt O_R(\widetilde{\pi}_k J_k(A)) \big\| \le \Big\| \meann_R \big| \widetilde{\pi}_k J_k(A) (\xi) \big|^2 d\xi \Big\|^\frac12 + \Big\| \meann_R \widetilde{\pi}_kJ_k(A) (\xi) d\xi \Big\| =: \alpha_R + \beta_R.$$ Splitting $R$ into a disjoint union of $\mathrm{K}$ small cubes $Q_j$ of the same measure, we find $$\alpha_R^2 \le \frac{1}{\mathrm{K}} \sum_{j=1}^\mathrm{K} \Big\| \meann_{Q_j} \big| \widetilde{\pi}_k J_k(A) (\xi) \big|^2 d\xi \Big\| \lesssim \sup_{Q \ \mathrm{small}} \big\| M\hskip-1pt O_Q(\widetilde{\pi}_kJ_k(A)) \big\|^2 + \frac{1}{\mathrm{K}} \sum_{j=1}^\mathrm{K} \beta_{Q_j}^2.$$ Moreover, we may assume that $$\mathrm{diam}(Q_j) \ge \frac12 4^k.$$ In particular, it suffices to prove that $\beta_Q \lesssim 2^{-k} \|J_k(A)\|_2$ for any such cube, also including $R$. Since $|p_k(x) - p_k(y)| \gtrsim 2^{-k}$ for any pair $(x,y)$ with $J_k(A)_{xy} \neq 0$, the desired inequality reduces to the simple fact that 
$$\Big| \meann_Q e^{2\pi i \langle \xi, p_k(x)-p_k(y) \rangle} d\xi \Big| \lesssim 2^{-k}$$ for every cube $Q$ with $\mathrm{diam}(Q) \gtrsim 4^k$. The proof for row $\mathrm{BMO}$-spaces is similar. \fin

Before proving Theorem \ref{Thm-Interpolation}, we just recall an auxiliary interpolation result which is not written explicitly in the literature, but might be known to experts. Consider the standard dyadic cubes in $\T^n$ and the corresponding martingale filtration in the von Neumann algebra $\A$. Construct $\mathrm{BMO}(\A)$, the martingale BMO space as defined by Pisier and Xu \cite{PX} equipped with its natural o.s.s. Let $\mathrm{BMO}_\circ(\A)$ be its mean-zero subspace, the image of the projection $J_\circ = \mathrm{id} - E_1$ $$J_\circ: \mathrm{BMO}(\A) \ni f \mapsto f - \int_{\T^n} f(x) dx \in \mathrm{BMO}_\circ(\A).$$ Recall $\displaystyle \|f\|_{\mathrm{BMO}(\A)} = \sup_{j \ge 1} \Big\| E_j \big| f - E_{j-1}(f) \big|^2 \Big\|_\infty^\frac12 \! \ge \max \Big\{ \big\| J_\circ(f) \big\|_{\mathrm{BMO}(\A)}, \|E_1(f)\|_\infty \Big\}$.

\begin{lemma} \label{Lemma-MeiMusat}
We have $$L_p^\circ(\A) \simeq_{\mathrm{cb}} \big[ \mathrm{BMO}_\A, L_2^\circ(\A) \big]_{\frac{2}{p}} \quad \mbox{where} \quad \mbox{$c_p \approx p$ as $p \to \infty$}.$$
\end{lemma}

\dem This is well known from \cite{Mei07,Mu}, the only novelty is the behavior of the constant $c_p$. As the inclusion $L_p^\circ(\A) \subset [\mathrm{BMO}_\A, L_2^\circ(\A)]_{2/p}$ is a complete contraction all we need to prove is that the interpolation space $[\mathrm{BMO}_\A, L_2^\circ(\A)]_{2/p}$ cb-embeds into $L_p^\circ(\A)$ with $c_p \approx p$. We may replace $\mathrm{BMO}_\A$ by the larger space $\mathrm{BMO}_\circ(\A)$ so that it suffices to cb-embed $[\mathrm{BMO}(\A), L_2(\A)]_{2/p}$ into $L_p(\A)$ and use $J_\circ$. By duality, we need the cb-embedding constant of $L_q(\A)$ into $[\mathrm{H}_1(\A), L_2(\A)]_{2/p}$ as $q \to 1$. Next, recalling that $L_2(\A) = \mathrm{H}_2^r(\A) = \mathrm{H}_2^c(\A)$ we deduce $$\mathrm{H}_q^r(\A) + \mathrm{H}_q^c(\A) \subset \big[ \mathrm{H}_1^r(\A), \mathrm{H}_2^r(\A) \big]_{\frac{2}{p}} + \big[ \mathrm{H}_1^c(\A), \mathrm{H}_2^c(\A) \big]_{\frac{2}{p}} \subset \big[ \mathrm{H}_1(\A), L_2(\A) \big]_{\frac{2}{p}}.$$ Moreover, these inclusions hold with absolute constants from \cite[Section 3]{Ra}. Thus our cb-embedding constant is controlled by the sharp growth of the constant in the noncommutative Burkholder-Gundy inequality $\|f\|_{\mathrm{H}_q} \le \alpha_q \|f\|_q$ \cite{PX}. We know from \cite{JRWZ} that $\alpha_q \approx 1/(q-1)$ as $q \to 1$. Therefore, we get $c_p \approx p$ as $p \to \infty$. \fin

\demInt The inclusion $$S_p(\R^n) = \big[ S_\infty(\R^n), S_2(\R^n) \big]_{\frac{2}{p}} \subset \big[ \mathrm{BMO}, S_2(\R^n) \big]_{\frac{2}{p}} =: \mathrm{X}_p(\R^n)$$ and its matrix $m$-amplification are clear. Therefore, according to Lemma \ref{Lemma-Interp-Density}, it suffices to prove the inequality below for any $A \in S_2(\R^n)$ and some constant $c_p \approx p$  
\begin{eqnarray} \label{Eq-Upper}
\|A\|_{S_p(\R^n)} \!\! & \le & \!\! c_p \hskip1pt \|A\|_{\mathrm{X}_p(\R^n)} \\ \nonumber \!\! & \approx & \!\! c_p \sup_{s \in \R} \max \Big\{ \big\| \Psi_{p,\varepsilon}(is) \big\|_{A_0}, \big\| \Psi_{p,\varepsilon}(1+is) \big\|_{A_1} \Big\}
\end{eqnarray}
with $(A_0,A_1) = (\mathrm{BMO}, S_2(\R^n))$. We only prove this inequality, since its matrix amplifications follow with the same proof and constant by \cite[Corollary 1.4]{PisAst}. Now note that $\|A - J_k(A)\|_2 \to 0$ as $k \to \infty$. Then, we use that $\pi_k \circ J_k(A)$ is mean-zero and the interpolation isomorphism in Lemma \ref{Lemma-MeiMusat} to deduce  
\begin{eqnarray*} 
\|A\|_{S_p(\R^n)} \!\! &= & \!\! \lim_{k \to \infty} \big\| J_k(A) \big\|_{S_p(\R^n)} \\ [3pt] \!\! &= & \!\! \lim_{k \to \infty} \big\| \pi_k \circ J_k(A) \big\|_{L_p^\circ(\A)} \lesssim \lim_{k \to \infty} \big\| \pi_k \circ J_k(A) \big\|_{[\mathrm{BMO}_\A, L_2^\circ(\A)]_{\frac2p}} \\ \!\! & \le & \!\! \lim_{k \to \infty} \sup_{s \in \R} \ \max \Big\{ \big\| \pi_k \circ J_k (\Psi_{p,\varepsilon}(is)) \big\|_{B_0}, \big\| \pi_k \circ J_k (\Psi_{p,\varepsilon}(1+is)) \big\|_{B_1} \Big\}
\end{eqnarray*}
with $(B_0,B_1) = (\mathrm{BMO}_\A, L_2^\circ(\A))$. The first inequality holds with $c_p \approx p$. The second one is possible since $\pi_k \circ J_k \circ \Psi_{p,\varepsilon} \in \mathcal{F}(B_0,B_1)$ and $\Psi_{p,\varepsilon}(2/p)=A$. More precisely, $\pi_k \circ J_k \circ \Psi_{p,\varepsilon}$ is clearly holomorphic on the strip and continuous on its closure. Moreover, Lemma \ref{Lemma-Interp-Density} ii) implies that its restriction to $\partial_0$ belongs to $C_0(\partial_0; \mathrm{BMO}_\A)$. Additionally, its restriction to $\partial_1$ belongs to $C_0(\partial_1; L_2^\circ(\A))$ since $\pi_k \circ J_k: S_2(\R^n) \to L_2^\circ(\A)$ is completely bounded. Altogether, this justifies that $\pi_k \circ J_k \circ \Psi_{p,\varepsilon} \in \mathcal{F}(B_0,B_1)$. Next, using Lemma \ref{Lemma-Interp-Density} ii) once more $$\sup_{s \in \R} \big\| J_k(\Psi_{p,\varepsilon}(is)) \big\|_2 \le \mathrm{M} \sup_{s \in \R} \big\| \Psi_{p,\varepsilon}(is) \big\|_{\mathrm{BMO}}$$ for some $\mathrm{M}>0$ independent of $k$. Then, Lemma \ref{Lemma-PikPi} gives for large enough $k$ 
\begin{eqnarray*}
\lefteqn{\hskip-20pt \sup_{s \in \R} \big\| \pi_k J_k(\Psi_{p,\varepsilon}(is)) \big\|_{B_0}} \\ \hskip20pt \!\!\! & \le & \!\!\!  \sup_{s \in \R} \big\| \pi J_k(\Psi_{p,\varepsilon}(is)) \big\|_{\mathrm{BMO}_\RR} + C_n 2^{-k} \mathrm{M} \sup_{s \in \R} \big\| \Psi_{p,\varepsilon}(is) \big\|_{\mathrm{BMO}} \\ \!\!\! & \le & \!\!\!  2 \sup_{s \in \R} \Big( \big\| \Psi_{p,\varepsilon}(is) \big\|_{A_0} + \big\| J_k^\perp(\Psi_{p,\varepsilon}(is)) \big\|_{S_2(\R^n)} \Big) \le 3 \sup_{s \in \R} \big\| \Psi_{p,\varepsilon}(is) \big\|_{A_0}
\end{eqnarray*}
with $J_k^\perp = \mathrm{id} - J_k$. The last inequality follows for large enough $k$ thanks to Lemma \ref{Lemma-Interp-Density} ii) and the dominated convergence theorem. On the other hand, the same inequality for $(A_1,B_1)$ is straightforward and the proof is complete. \fin 

\begin{remark} \label{Rem-RowColInterpolation}
\emph{Theorem \ref{Thm-Interpolation} easily extends to $$[\mathrm{BMO}, S_p(\R^n)]_{\frac{p}{q}} \simeq_{\mathrm{cb}} S_q(\R^n) \, \simeq_{\mathrm{cb}} [\mathrm{BMO}_r, S_2(\R^n)]_{\frac2q} \cap [\mathrm{BMO}_c, S_2(\R^n)]_{\frac2q}$$ for $2 \le p < q < \infty$ by Wolff's interpolation as in \cite{Mu}. The second isomorphism follows as in Theorem \ref{Thm-Interpolation} using $(A_0,B_0) = (\mathrm{BMO}_\dag, \mathrm{BMO}_{\A}^\dag)$ for $\dag = r,c$ separately.}
\end{remark}

\begin{remark} \label{Rmk-Interp-Semigroups}
\emph{An equivalent BMO norm on $\mathcal{B}(L_2(\R^n))$ can be constructed by means of a Markov semigroup which arises by transferring the $n$-dimensional heat semigroup. More precisely, one could set $$\S_t(A) = \Big( \exp(-t|x-y|^2) A_{xy} \Big) \quad \mbox{and} \quad \|A\|_{\mathrm{BMO}_\S^c} = \sup_{t > 0} \Big\| \Big( \S_t (|A|^2) - |\S_t(A)|^2 \Big)^\frac12 \Big\|_\infty.$$ The techniques in \cite{JM} show as well that such BMO interpolates with the $S_p$-scale. Our definition though is more natural for matrix algebras and the interpolation argument is elementary. Moreover, it generalizes for anisotropic forms of our BMO space which is crucial in \cite{CGPT}, whereas no Markov semigroup seems to work in that setting. Other BMO spaces were considered in \cite{CJSZ} for divided differences.}
\end{remark}

\section{\bf \large H\"ormander-Mikhlin-Schur multipliers}

In this section we prove our main result. In fact, we shall prove a stronger statement which only assumes a fractional regularity order $\sigma = \frac{n}{2} + \varepsilon$ for any fixed $\varepsilon > 0$. However, for the sake of clarity we first prove Theorem A and explain the modifications afterwards. The proof will be divided into several steps.  

\noindent \textbf{1) Reduction to BMO.} Since $S_M^* = S_{\overline{M}}$ and our HMS conditions in Theorem A are stable by complex conjugation, we may assume with no loss of generality that $p \ge 2$. The boundedness for $p=2$ clearly holds since the HMS condition readily implies that $M$ is essentially bounded on $\R^n \hskip-2pt \times \R^n$. It suffices to prove the inequality 
\begin{equation} \label{Eq-BMOEndpoint}
\big\| S_M: S_\infty(\R^n) \to \mathrm{BMO} \big\| \ \le \ C ||| M |||_{\mathrm{HMS}}.
\end{equation}
Theorem A follows from \eqref{Eq-BMOEndpoint} by Theorem \ref{Thm-Interpolation} and Riesz-Thorin interpolation. 

\noindent \textbf{2) Twisted Fourier multipliers.} Define 
\begin{equation} \label{Eq-McMr}
M_r(x,y) := M(y-x,y) \quad \mbox{and} \quad M_c(x,y) := M(x,x-y).
\end{equation}
Then we construct the following linear maps on $L_2(\RR)$ 
\begin{eqnarray*}
\widetilde{T}_{M_r}(f)(z) & := & \Big( \int_{\R^n} M_r(\xi,y) \widehat{f}_{xy}(\xi) e^{2\pi i \langle z, \xi \rangle}\, d\xi \Big), \\
\widetilde{T}_{M_c}(f)(z) & := & \Big( \int_{\R^n} M_c(x,\xi) \widehat{f}_{xy}(\xi) e^{2\pi i \langle z, \xi \rangle}\, d\xi \Big).
\end{eqnarray*}
The use of both maps (not just one) is critical. Similar maps were already considered in \cite{PRS,PRo} and denominated twisted Fourier multipliers. They act on the entry $f_{xy}$ as the Fourier multiplier with symbol $M_r(\cdot, y)$ or $M_c(x,\cdot)$ respectively. Note that the finiteness of $|||M|||_{\mathrm{HMS}}$ in the statement readily implies that $M_r(\cdot, y)$ and $M_c(x,\cdot)$ are Mikhlin multipliers and therefore well-defined on $L_\infty(\R^n)$ to $\mathrm{BMO}(\R^n)$. More than that, we shall justify in points 3), 4) and 5) below that the corresponding twisted multipliers are well-defined on $L_\infty(\RR)$ with values in row/column $\mathrm{BMO}_\RR$. In this part of the proof, we shall concentrate on the inequalities  
\begin{eqnarray} 
\label{Eq-Transf1}
\big\| S_M: S_\infty(\R^n) \to \mathrm{BMO}_r \big\| \!\! & \le & \!\! \big\| \widetilde{T}_{M_r}: L_\infty(\RR) \to \mathrm{BMO}_{\RR}^r \big\|, \\  
\label{Eq-Transf2}
\big\| S_M: S_\infty(\R^n) \to \mathrm{BMO}_c \big\| \!\! & \le & \!\! \big\| \widetilde{T}_{M_c}: L_\infty(\RR) \to \mathrm{BMO}_{\RR}^c \big\|.
\end{eqnarray}
Let us first focus on the second inequality \eqref{Eq-Transf2}. Using the $*$-homomorphism $\pi$ above, we shall use an additional intertwining identity between the Schur multiplier $S_M$ and the twisted Fourier multiplier associated to $M_c$
\begin{eqnarray} \label{Eq-Intertwining2}
\hskip20pt \widetilde{T}_{M_c}(\pi(A))(z) \!\! & = & \!\! \Big( T_{M_c(x,\cdot)} \big( \exp(2 \pi i \langle x-y, \cdot \rangle) \big)(z) A_{xy} \Big) \\ \nonumber \!\! & = & \!\! \Big( M_c(x,x-y) \exp(2 \pi i \langle x-y, z \rangle) A_{xy} \Big) \ = \ \pi (S_M(A)).
\end{eqnarray}
Then, we get 
\begin{eqnarray*}
\|S_M(A)\|_{\mathrm{BMO}_c} \!\! & = & \!\! \big\| \pi (S_M(A)) \big\|_{\mathrm{BMO}_\RR^c} \\ \!\! & = & \!\! \big\| \widetilde{T}_{M_c}(\pi(A)) \big\|_{\mathrm{BMO}_\RR^c} \ \le \ \big\| \widetilde{T}_{M_c}: L_\infty(\RR) \to \mathrm{BMO}_{\RR}^c \big\| \, \|\pi(A)\|_\infty.
\end{eqnarray*}
Since $\pi$ acts isometrically on $S_\infty(\R^n)$, this proves inequality \eqref{Eq-Transf2}. The proof of \eqref{Eq-Transf1} is essentially the same. The only change consists in replacing $\pi$ by its adjoint
$$\sigma (A)(\xi) := \Big( \exp \big( 2\pi i \langle \xi, y-x \rangle \big) A_{xy} \Big).$$ 
Indeed, the intertwining identity \eqref{Eq-Intertwining2} still holds for $(M_r,\sigma)$ in place of $(M_c,\pi)$.

\noindent \textbf{3) Calder\'on-Zygmund kernels.} We claim that
\begin{eqnarray} 
\label{Eq-HM1}
\big\| \widetilde{T}_{M_r}: L_\infty(\RR) \to \mathrm{BMO}_{\RR}^r \big\| \!\! & \lesssim & \!\! \sum_{|\gamma| \le [\frac{n}{2}] +1} \Big\| |x-y|^{|\gamma|} \big| \partial_x^\gamma M(x,y) \big| \Big\|_\infty, \\ [-3pt]
\label{Eq-HM2}
\big\| \widetilde{T}_{M_c}: L_\infty(\RR) \to \mathrm{BMO}_{\RR}^c \big\| \!\! & \lesssim & \!\! \sum_{|\gamma| \le [\frac{n}{2}] +1} \Big\| |x-y|^{|\gamma|} \big| \partial_y^\gamma M(x,y) \big| \Big\|_\infty.
\end{eqnarray}
This implies Theorem A. The proof requires us to refine well-known arguments from classical Calder\'on-Zygmund theory in the noncommutative realm. Again, we prove first the column endpoint. We formally have
\begin{eqnarray*}
\widetilde{T}_{M_c}(f)(z) \!\! & = & \!\! \Big( \int_{\R^n} M_c(x,\xi) \widehat{f}_{xy}(\xi) e^{2\pi i \langle z, \xi \rangle}\, d\xi \Big) \\ [7pt]
\!\! & = & \!\! \Big( \int_{\R^n} \underbrace{\big[ M_c(x,\cdot) \big]^{\vee}(z-w)}_{k_c(x,z-w)} f_{xy}(w) \, dw \Big) \\ \!\! & = & \!\!  \Big( k_c(x,\cdot) * f_{xy} \Big)(z) \, = \, \int_{\R^n} K_c(z-w) \cdot f(w) \, dw
\end{eqnarray*}    
with $K_c: \R^n \setminus \{0\} \to \mathcal{B}(L_2(\R^n))$ the \lq\lq diagonal-valued" function $$(K_c(z)h)(x) = k_c(x,z) h(x).$$ Although it is not meaningful, we shall abusively use the matrix-notation $$K_c = \mathrm{diag} \Big( k_c(x,\cdot) \Big).$$
The kernel $K_c$ should be understood as an operator-valued distribution which agrees with a locally integrable operator-valued function on $\R^n \setminus \{0\}$. A careful analysis of these kernels is analogous to that from \cite[Section 2]{GJP2}, where CZ kernels on quantum Euclidean spaces were investigated. In the present context, it helps to understand the behavior of the \lq\lq diagonal entries" $k_c(x,\cdot) = M_c(x,\cdot)^{\vee}$ of the kernel $K_c$. Since we allow $M$ to be singular at $x=y$ it turns out that $M_c(x,\cdot)$ may be singular at $0$. This means that its inverse Fourier transform must be understood as a distribution which is locally integrable outside the origin. In particular, the kernel representation above is meaningful when $z \notin \mathrm{supp}_{\R^n} f$, the Euclidean support of the matrix-valued function $f$. Under this assumption (standard in CZ theory) the kernel representation above does not require to deal with distributions taking values in operator algebras. More precisely, we shall only use this kernel representation for $z \in \R^n \setminus \mathrm{supp}_{\R^n} f$, which is exactly the framework of the noncommutative CZ theory \cite{CCP,JMP1,Pa1}. In our use of CZ methods, it will be crucial to weaken the kernel condition from  \cite[Lemma 2.3]{JMP1}. As usual in this context, we begin by noticing that $$\big\| \widetilde{T}_{M_c} f \big\|_{\mathrm{BMO}_\RR^c} \sim \sup_{Q \in \Q} \inf_{\alpha_Q \in \B(L_2(\R^n))} \Big\| \Big( \frac{1}{|Q|} \int_Q \big| \widetilde{T}_{M_c} f(z) - \alpha_Q \big|^2 dz \Big)^{\frac12} \Big\|_{\B(L_2(\R^n))},$$ where $\Q$ is the set of Euclidean balls as above. Next, fix $Q \in \Q$ nearly attaining the above sup and decompose $f = f_1 + f_2$ with $f_1 = f \chi_{5Q}$ and $f_2 = f - f_1$. Picking $$\alpha_Q = \meann_Q \widetilde{T}_{M_c}(f_2)(w) dw,$$ we obtain the following estimate 
\begin{eqnarray*}
\big\| \widetilde{T}_{M_c} f \big\|_{\mathrm{BMO}_\RR^c} \!\! & \lesssim & \!\! \Big\| \Big( \mean_Q \big| \widetilde{T}_{M_c} f_1(z) \big|^2 dz \Big)^{\frac12} \Big\|_\infty \\ \!\! & + & \!\! \Big\| \Big( \mean_Q \Big| \mean_Q \widetilde{T}_{M_c} f_2(z) - \widetilde{T}_{M_c} f_2(w) \, dw \Big|^2 dz \Big)^{\frac12} \Big\|_\infty \ = \ \mathrm{A}+\mathrm{B}.
\end{eqnarray*}
On the other hand, we recall the following inequalities:
\begin{eqnarray}
\label{Eq-L2Bdness}
\big\| \widetilde{T}_{M_c}f \big\|_{L_2(\RR)} \!\! & \le & \!\! \|M_c\|_\infty \|f\|_{L_2(\RR)}, \\
\label{Eq-LinftyL2Col}
\hskip30pt \Big\| \int_{\R^n} \big| \widetilde{T}_{M_c}f(z) \big|^2 \, dz \Big\|_{\B(L_2(\R^n))}^{\frac12} \!\! & \le & \!\! \|M_c\|_\infty \Big\| \int_{\R^n} \big| f(z) \big|^2 \, dz \Big\|_{\B(L_2(\R^n))}^{\frac12}.   
\end{eqnarray}
Inequality \eqref{Eq-L2Bdness} is straightforward from the Hilbert-valued Plancherel theorem and \eqref{Eq-LinftyL2Col} follows from it, since $M_c$-twisted multipliers are right module maps and thus $$\Big\| \int_{\R^n} \hskip-2pt \big| \widetilde{T}_{M_c}f(z) \big|^2 dz \Big\|_\infty \hskip-2pt  = \hskip-2pt \sup_{\|E\|_{S_2}=1} \int_{\R^n} \hskip-2pt  \mathrm{tr} \Big( \big| \widetilde{T}_{M_c}[fE](z) \big|^2 \Big) dz = \hskip-2pt \sup_{\|E\|_{S_2}=1} \big\| \widetilde{T}_{M_c}[fE] \big\|_2^2,$$ from where our claim follows. By \eqref{Eq-LinftyL2Col} 
$$\Big\| \Big( \mean_Q \big| \widetilde{T}_{M_c} f_1(z) \big|^2 dz \Big)^{\frac12} \Big\|_\infty \le \frac{\|M_c\|_\infty}{\sqrt{|Q|}} \Big\| \Big( \int_{\R^n} \big| f_1(z) \big|^2 dz \Big)^{\frac12} \Big\|_\infty.$$ Since $f_1 = f \chi_{5Q}$ we immediately deduce
\begin{equation} \label{Eq-A}
\mathrm{A} \, \le \, 5^{\frac{n}{2}} \|M_c\|_\infty \|f\|_\infty.
\end{equation}
To estimate B, the operator-Jensen inequality and the kernel representation yield
\begin{eqnarray*}
\mathrm{B} \!\! & \le & \!\! \sup_{Q \in \Q} \Big\| \Big( \mean_{Q \times Q} \big| \widetilde{T}_{M_c} f_2(z) - \widetilde{T}_{M_c} f_2(w) \big|^2 dwdz \Big)^{\frac12} \Big\|_\infty \\ 
\!\! & \le & \!\! \sup_{Q \in \Q} \sup_{w,z \in Q} \Big\| \int_{\R^n \setminus 5Q} \big( K_c(z-s) - K_c(w-s) \big) f(s)  \, ds \Big\|_\infty.
\end{eqnarray*}
It is very tempting to extract the $L_\infty$-norm of $f$ by putting the norm inside the integral. Since $|w-s| > 2 |z-w|$ for $s,w,z$ as above, we could use $s \mapsto w-s$ and $u = w-z$ to put the region of integration inside the set $|u| < \frac12 |s|$. This approach from \cite[Lemma 2.3]{JMP1} is however too restrictive for our CZ kernel under the given (sharp) HM conditions. Instead, we set $f_w(s) = f(w-s)$ and $\Gamma = Q - Q$ to obtain the following inequality
\begin{equation} \label{Eq-Hormander}
\mathrm{B} \, \le \, \sup_{(u,w) \in \Gamma \times Q} \underbrace{\Big\| \int_{\R^n \setminus (w - 5Q)} \big( K_c(s-u) - K_c(s) \big) f_w(s)  \, ds \Big\|_\infty}_{\mathrm{B}_{(u,w)}}. 
\end{equation}
Then, we use a Littlewood-Paley decomposition of $K_c$ before loosing cancellation.

\noindent \textbf{4) H\"ormander-Mikhlin estimates.}  Let us continue with the proof of inequality  \eqref{Eq-HM2}. According to \eqref{Eq-A}, we have $\mathrm{A}/\|f\|_\infty \lesssim \mathrm{RHS} \hskip1pt \eqref{Eq-HM2}$. It therefore remains to control the term $\mathrm{B}$. Let $\phi: \R^n \to \R_+$ be a smooth function which is identically one in the unit ball and vanishes for $|\xi| > 2$. The functions $\psi_j(\xi) = \phi(2^{-j}\xi) - \phi(2^{1-j}\xi)$ define a Littlewood-Paley partition of unity. That is to say 
\begin{equation} \label{Eq-LP}
\mathrm{supp} (\psi_j) \subset \big\{2^{j-1} \le |\xi| \le 2^{j+1}\big\} \quad \mbox{and} \quad \sum_{j \in \Z} \psi_j(\xi) = 1 \quad \mbox{for $\xi \neq 0$}.
\end{equation}
Decompose $$K_c = \sum_{j \in \Z} K_{cj}$$ with the pieces given by $$K_{cj} = \mathrm{diag} \Big( k_{cj}(x,\cdot) \Big) \quad \mbox{and} \quad k_{cj}(x,z) = \big[\psi_j M_c(x,\cdot) \big]^\vee(z).$$ 
Let us recall that Plancherel theorem gives 
\begin{equation} \label{Eq-Plancherel} 
\int_{\R^n} \big| (-2 \pi i z)^\gamma k_{cj}(x,z) \big|^2 dz = \int_{\R^n} \big| \partial_\xi^\gamma \big( M_{c}(x,\xi) \psi_j(\xi) \big) \big|^2 d\xi
\end{equation}
for any multi-index $\gamma$. Given a fixed $(u,w) \in \Gamma \times Q$, we have 
\begin{eqnarray*}
\mathrm{B}_{(u,w)} \!\!\!\! & \le & \!\!\!\!\!\! \sum_{\begin{subarray}{c} j \in \Z \\ 2^j |u| > 1\end{subarray}} \Big\| \int_{\R^n \setminus (w - 5Q)} \big( K_{cj}(s-u) - K_{cj}(s) \big) f_w(s)  \, ds \Big\|_\infty \\ [-3pt] \!\!\!\! & + & \!\!\!\!\!\! \sum_{\begin{subarray}{c} j \in \Z \\ 2^j |u| \le 1 \end{subarray}} \Big\| \int_{\R^n \setminus (w - 5Q)} \big( K_{cj}(s-u) - K_{cj}(s) \big) f_w(s)  \, ds \Big\|_\infty =: \mathrm{B}_{(u,w)}^1 + \mathrm{B}_{(u,w)}^2.
\end{eqnarray*}
By the triangle inequality
\begin{eqnarray*}
\mathrm{B}_{(u,w)}^1 \!\!\!\! & \le & \!\!\!\!\!\! \sum_{\begin{subarray}{c} j \in \Z \\ 2^j |u| > 1\end{subarray}} \Big\| \int_{\R^n \setminus (w - 5Q)} K_{cj}(s-u) f_w(s) \, ds \Big\|_\infty \\ [-3pt] \!\!\!\! & + & \!\!\!\!\!\! \sum_{\begin{subarray}{c} j \in \Z \\ 2^j |u| > 1\end{subarray}} \Big\| \int_{\R^n \setminus (w - 5Q)} K_{cj}(s) f_w(s) \, ds \Big\|_\infty =:  \sum_{\begin{subarray}{c} j \in \Z \\ 2^j |u| > 1\end{subarray}} \mathrm{B}_{(u,w)}^{11}(j) + \mathrm{B}_{(u,w)}^{12}(j).
\end{eqnarray*}
Next, letting $\sigma = [\frac{n}{2}]+1$ and using H\"older inequality 
$$\mathrm{B}_{(u,w)}^{11}(j) \, \le \, \Big\| \int_{\R^n} |s-u|^{2\sigma} |K_{cj}^*(s-u)|^2 \, ds \Big\|_\infty^{\frac12} \Big\| \int_{|u| < \frac12|s|} |s-u|^{-2\sigma} |f_w(s)|^2 \, ds \Big\|_\infty^{\frac12}.$$
Since $|s+u| \ge |s| - |u| > |u|$ when $|u| < \frac12|s|$ and $K_{cj} K_{cj}^* = K_{cj}^* K_{cj}$
\begin{equation} \label{Eq-B11}
\mathrm{B}_{(u,w)}^{11}(j) \, \le \, \Big\| \int_{\R^n} |s|^{2\sigma} |K_{cj}(s)|^2 \, ds \Big\|_\infty^{\frac12} \Big( \int_{|u| < |s|} |s|^{-2\sigma} \, ds \Big)^{\frac12} \|f\|_\infty. 
\end{equation}
Using $\| m \|_{\mathrm{HM}} := \sum_{|\gamma| \le \sigma} \big\| |\xi|^{|\gamma|} \partial_\xi^\gamma m(\xi) \big\|_\infty$ and \eqref{Eq-Plancherel}, we find the estimate 
\begin{eqnarray*}
\lefteqn{\hskip-20pt \Big\| \int_{\R^n} |s|^{2\sigma} |K_{cj}(s)|^2 \, ds \Big\|_\infty} \\ \!\!\! & \lesssim & \!\!\! \sup_{x \in \R^n} \ \sup_{|\gamma| = \sigma} \int_{\R^n} \big| \partial_\xi^\gamma \big( \psi_j(\xi) M_{c}(x,\xi) \big) \big|^2 d\xi \\ \!\!\! & \lesssim & \!\!\! \sup_{x \in \R^n} \sum_{|\alpha| + |\beta| = \sigma} \int_{\R^n} \big| \partial_\xi^\alpha  \psi_j(\xi) \partial_\xi^\beta M_{c}(x,\xi) \big|^2 d\xi \\ \!\!\! & \lesssim & \!\!\! \sup_{x \in \R^n} \sum_{|\alpha| + |\beta| = \sigma} 2^{-2j |\alpha|} \int_{\R^n} \big| \partial_\xi^\alpha \psi(2^{-j}\xi) |\xi|^{-|\beta|} \big|^2 d\xi \ \big\| M_c(x,\cdot) \big\|_{\mathrm{HM}}^2.
\end{eqnarray*}
Since $\mathrm{supp} (\psi) \subset \{\frac12 \le |\xi| \le 2\}$, this is bounded by $2^{j(n-2\sigma)} \sup_x \| M_c(x,\cdot) \|_{\mathrm{HM}}^2$. Thus, by \eqref{Eq-B11} we obtain
\begin{eqnarray*}
\lefteqn{\hskip-20pt \sum_{\begin{subarray}{c} j \in \Z \\ 2^j |u| > 1\end{subarray}} \mathrm{B}_{(u,w)}^{11}(j) \, \lesssim \, |u|^{\frac{n}{2}-\sigma} \sum_{\begin{subarray}{c} j \in \Z \\ 2^j |u| > 1\end{subarray}} 2^{j (\frac{n}{2} - \sigma)} \sup_{x \in \R^n} \big\| M_c(x,\cdot) \big\|_{\mathrm{HM}} \, \|f\|_\infty} \\ \hskip20pt \!\!\!\! & \lesssim & \!\!\!\! \sup_{x \in \R^n} \big\| M_c(x,\cdot) \big\|_{\mathrm{HM}} \, \|f\|_\infty = \!\! \sum_{|\gamma| \le [\frac{n}{2}] +1} \Big\| |x-y|^{|\gamma|} \big| \partial_y^\gamma M(x,y) \big| \Big\|_\infty \|f\|_\infty
\end{eqnarray*}
since $\partial_\xi^\gamma M_c(x,\xi) = (-1)^{|\gamma|} \partial_\xi^\gamma M(x,x-\xi)$. The sum of $\mathrm{B}_{(u,w)}^{12}(j)$'s admits the same bound with the same argument. Altogether, we have proved the following inequality 
\begin{equation} \label{Eq-B1}
\mathrm{B}_{(u,w)}^1 \, \lesssim \, \sum_{|\gamma| \le [\frac{n}{2}] +1} \Big\| |x-y|^{|\gamma|} \big| \partial_y^\gamma M(x,y) \big| \Big\|_\infty \|f\|_\infty. 
\end{equation}
To estimate $\mathrm{B}_{(u,w)}^2$ we set $A_w(Q) = \R^n \setminus (w - 5Q)$
\begin{eqnarray*}
\mathrm{B}_{(u,w)}^2 \!\!\!\! & \le & \!\!\!\!\!\! \sum_{\begin{subarray}{c} j \in \Z \\ 2^j |u| \le 1\end{subarray}} \Big\| \int_{A_w(Q), |s| < 2^{-j}} \big( K_{cj}(s-u) - K_{cj}(s) \big) f_w(s) \, ds \Big\|_\infty \\ [-3pt]
\!\!\!\! & + & \!\!\!\!\!\! \sum_{\begin{subarray}{c} j \in \Z \\ 2^j |u| \le 1\end{subarray}} \Big\| \int_{A_w(Q), |s| \ge 2^{-j}} \big( K_{cj}(s-u) - K_{cj}(s) \big) f_w(s) \, ds \Big\|_\infty \!\! = \mathrm{B}_{(u,w)}^{21} + \mathrm{B}_{(u,w)}^{22}.
\end{eqnarray*}
The $j$-th term of $\mathrm{B}_{(u,w)}^{21}$ is dominated by 
\begin{eqnarray*}
\mathrm{B}_{(u,w)}^{21}(j) \!\! & \le & \!\! \Big\| \int_{|u| < \frac12 |s|} \big| K_{cj}(s-u)^* - K_{cj}(s)^* \big|^2 \, ds \Big\|_\infty^{\frac12} \Big\| \int_{|s| < 2^{-j}} |f_w(s)|^2 \, ds \Big\|_\infty^{\frac12} \\ \!\! & \lesssim & \!\! 2^{-\frac{jn}{2}} \sup_{x \in \R^n} \Big( \int_{|u| < \frac12 |s|} \big| k_{cj}(x,s-u) - k_{cj}(x,s) \big|^2 \, ds \Big)^{\frac12} \|f\|_\infty \\ \!\! & = & \!\! 2^{-\frac{jn}{2}} \sup_{x \in \R^n} \Big( \int_{|u| < \frac12 |s|} \Big| \int_0^1 \frac{d}{dr} k_{cj}(x,s-ru) \, dr \Big|^2 \, ds \Big)^{\frac12} \|f\|_\infty \\ \!\! & \lesssim & \!\! 2^{-\frac{jn}{2}} \sup_{x \in \R^n} \Big( \int_{\R^n} \big| \nabla k_{cj}(x,s) \big|^2 \, ds \Big)^{\frac12} |u| \, \|f\|_\infty \ \lesssim \ 2^j |u| \, \|M_c\|_\infty \|f\|_\infty 
\end{eqnarray*}
since the Plancherel theorem and $\mathrm{supp} (\psi) \subset \{\frac12 \le |\xi| \le 2\}$ give $$\int_{\R^n} \big| \nabla k_{cj}(x,s) \big|^2 \, ds \, = \, 4 \pi^2 \int_{\R^n} |\xi|^2 \big| M_c(x,\xi) \psi_j(\xi) \big|^2 \, d\xi \, \lesssim \, 2^{j(n+2)} \|M_c\|_\infty^2.$$ In particular, we get the following bound for the sum of $\mathrm{B}_{(u,w)}^{21}(j)$'s
\begin{equation} \label{Eq-B21}
\sum_{2^j|u| \le 1} \mathrm{B}_{(u,w)}^{21}(j) \lesssim \|M_c\|_{\infty} \|f\|_\infty. 
\end{equation}
As for $\mathrm{B}_{(u,w)}^{22}(j)$ we proceed similarly using the weights $|s|^{\pm \sigma}$ to obtain
$$\mathrm{B}_{(u,w)}^{22}(j) \lesssim |u| \sup_{x \in \R^n} \Big( \underbrace{\int_{\R^n} |s|^{2\sigma} \big| \nabla k_{cj}(x,s) \big|^2 \, ds}_{J_j(x)} \Big)^{\frac12} \Big\| \int_{|s| \ge 2^{-j}} |s|^{-2\sigma} |f_w(s)|^2 \, ds \Big\|_\infty^{\frac12}.$$
By Plancherel, the HM condition for $M_c$ and $\mathrm{supp} (\psi) \subset \{\frac12 \le |\xi| \le 2\}$, we get
\begin{eqnarray*}
J_j(x) \!\!\!\! & \lesssim & \!\! \sum_{|\gamma| = \sigma} \sum_{k=1}^n \int_{\R^n} \big| \partial_\xi^\gamma [ \xi_k \psi_j(\xi) M_c(x,\xi) ] \big|^2 d\xi \\ \!\!\!\! & \lesssim & \!\! \sum_{|\gamma| = \sigma} \int_{\R^n} |\xi|^2 \big| \partial_\xi^\gamma [ \psi_j(\xi) M_c(x,\xi) ] \big|^2 d\xi \\ \!\!\!\! & + & \!\!\!\! \sum_{|\gamma| = \sigma-1} \int_{\R^n} \big| \partial_\xi^\gamma [ \psi_j(\xi) M_c(x,\xi) ] \big|^2 d\xi \ \lesssim \ 2^{j(n+2-2\sigma)} \big\|M_c(x,\cdot) \big\|_{\mathrm{HM}}^2.   
\end{eqnarray*}
Since $2\sigma > n$, we obtain the following bound for the sum of $\mathrm{B}_{(u,w)}^{22}(j)$'s
\begin{equation} \label{Eq-B22}
\sum_{2^j|u| \le 1} \mathrm{B}_{(u,w)}^{22}(j) \lesssim \sup_{x \in \R^n} \big\| M_c(x,\cdot) \big\|_{\mathrm{HM}} \|f\|_\infty. 
\end{equation}
Combining \eqref{Eq-B21} and \eqref{Eq-B22}, we conclude 
\begin{equation} \label{Eq-B2}
\mathrm{B}_{(u,w)}^2 \, \lesssim \, \sum_{|\gamma| \le [\frac{n}{2}] +1} \Big\| |x-y|^{|\gamma|} \big| \partial_y^\gamma M(x,y) \big| \Big\|_\infty \|f\|_\infty.
\end{equation}

\noindent \textbf{5) Conclusion.} According to \eqref{Eq-A}, \eqref{Eq-Hormander}, \eqref{Eq-B1} and \eqref{Eq-B2}, we get
$$\big\| \widetilde{T}_{M_c}: L_\infty(\RR) \to \mathrm{BMO}_{\RR}^c \big\| \, \lesssim \, \sum_{|\gamma| \le [\frac{n}{2}] +1}^{\null} \Big\| |x-y|^{|\gamma|} \big| \partial_y^\gamma M(x,y) \big| \Big\|_\infty.$$ This proves \eqref{Eq-HM2} and \eqref{Eq-HM1} is analogous. Namely, using $(K_r(z)h)(y) =  k_r(z,y) h(y)$, it now acts by right multiplication. In particular, \eqref{Eq-LinftyL2Col} should be replaced by its row analog, which is used to get the same estimate as in \eqref{Eq-A}. The rest of the proof of \eqref{Eq-HM1} is essentially the same as that of \eqref{Eq-HM2}, with the only difference being that we get the upper bound $$\sup_{y \in \R^n} \big\| M_r(\cdot, y) \big\|_{\mathrm{HM}} \, = \, \sum_{|\gamma| \le [\frac{n}{2}] + 1} \Big\| |x-y|^{|\gamma|} \partial_x^\gamma M(x,y) \Big\|_\infty.$$ Once we get \eqref{Eq-HM1} \hskip-2pt + \hskip-2pt \eqref{Eq-HM2}, the proof follows from \eqref{Eq-Transf1} \hskip-2pt + \hskip-2pt \eqref{Eq-Transf2} and Theorem \ref{Thm-Interpolation}. \fin    

\begin{remark} \label{Rmk-ThmA-cb} 
\emph{By property i) in Section \ref{Sect-Schur}, Theorem A holds as well in the category of operator spaces. Moreover, all our $L_\infty \to \mathrm{BMO}$ estimates hold in the cb-setting with the given o.s.s. in the different BMO spaces. The proofs follow as above by taking matrix amplifications at each step of the original argument. This is relevant in \cite{CGPT}, where we work with discrete groups and \cite[Theorem 1.18]{LdlS} does not apply.}
\end{remark}

Consider now the Littlewood-Paley partition of unity $\{\psi_j\}_{j \in \Z}$ from \eqref{Eq-LP} and set $\psi(\xi) = \phi(\xi) - \phi(2\xi)$, so that $\psi_j(\xi) = \psi(2^{-j} \xi)$. Given a symbol $M: \R^n \! \times \R^n \to \C$ we construct the functions $M_r$ and $M_c$ as in \eqref{Eq-McMr}. Let $W_{q\sigma}$ denote the inhomogeneous $L_q$-Sobolev space in $\R^n$ of order $\sigma$ 
$\| f \|_{W_{q\sigma}} := \| ( (1 + |\cdot| )^{\sigma} f^\wedge )^\vee \|_{L_q(\R^n)},$
and define the Sobolev $\mathrm{HMS}_{q \sigma}$ norm as follows 
\begin{eqnarray} \label{Eq-SobolevNorm}
\hskip30pt \big\bracevert \hskip-2pt M \hskip-2pt \big\bracevert_{\hskip-3pt \mathrm{HMS}_{q \sigma}} \hskip-2pt \!\!\!\! & := & \!\!\!\!\! \hskip-2pt  \sup_{\begin{subarray}{c} j \in \Z \\ x,y \in \R^n \end{subarray}} \big\| \psi M_{r}(2^j \cdot,y) \big\|_{W_{q \sigma}} \hskip-5pt + \big\| \psi M_c(x,2^j \cdot) \big\|_{W_{q \sigma}} \\ \nonumber \hskip-2pt \!\!\!\! & = & \!\!\!\!\! \hskip-2pt \sup_{\begin{subarray}{c} j \in \Z \\ x,y \in \R^n \end{subarray}} \big\| \underbrace{\psi(\cdot - y)  M(2^j \cdot,2^j y)}_{\Psi_{M,j}(\cdot,y)} \big\|_{W_{q \sigma}} \hskip-5pt + \big\| \underbrace{\psi(x - \cdot) M(2^j x,2^j \cdot)}_{\Psi_{M,j}(x,\cdot)} \big\|_{W_{q \sigma}}.
\end{eqnarray}
It is a simple exercise to check that 
\begin{eqnarray} \label{Eq-HMSComparison}
\|M\|_\infty \le \big\bracevert \hskip-2pt M \hskip-2pt \big\bracevert_{\hskip-3pt \mathrm{HMS}_{2 \sigma}} \le ||| M |||_{\mathrm{HMS}} \quad \mbox{for any $\frac{n}{2} < \sigma < [\frac{n}{2}] + 1$}.
\end{eqnarray} 

\begin{TheoB}[Sobolev form] 
If $1 < p < \infty$ and $\sigma = \frac{n}{2}+\varepsilon$
$$\big\| S_M \hskip-2pt: S_p(\R^n) \to S_p(\R^n) \big\|_{\mathrm{cb}} \, \le \, C_\varepsilon \frac{p^2}{p-1} \big\bracevert \hskip-2pt M \hskip-2pt \big\bracevert_{\hskip-3pt \mathrm{HMS}_{2 \sigma}}.$$
\end{TheoB} 

\ske According to the first inequality in \eqref{Eq-HMSComparison}, points 1) and 2) in the proof of Theorem A remain the same in this case, just replacing the HMS norm in \eqref{Eq-BMOEndpoint} by the $\mathrm{HMS}_{2 \sigma}$ analog in the statement. Therefore, all we need to do is to improve our estimates in \eqref{Eq-HM1} and \eqref{Eq-HM2} as follows
\begin{eqnarray} 
\label{Eq-HM3}
\big\| \widetilde{T}_{M_r}: L_\infty(\RR) \to \mathrm{BMO}_{\RR}^r \big\|_{\mathrm{cb}} \!\! & \lesssim & \!\! \sup_{(j,y) \in \Z \times \R^n} \big\| \psi M_{r}(2^j \cdot,y) \big\|_{W_{2\sigma}}, \\ 
\label{Eq-HM4}
\big\| \widetilde{T}_{M_c}: L_\infty(\RR) \to \mathrm{BMO}_{\RR}^c \big\|_{\mathrm{cb}} \!\! & \lesssim & \!\! \sup_{(j,x) \in \Z \times \R^n} \big\| \psi M_c(x,2^j \cdot) \big\|_{W_{2\sigma}}.
\end{eqnarray}
To that end we argue as in point 3) and part of 4). Given a fixed $(u,w) \in \Gamma \times Q$ we are reduced to estimating $\mathrm{B}_1 + \mathrm{B}_2$. For simplicity, we no longer specify the dependence of these terms on the pair $(u,w)$. If $\sigma = \frac{n}{2} + \varepsilon$ we note that 
\begin{eqnarray*}
\mathrm{B}_{11}(j) \!\! & \le & \!\! \Big\| \int_{\R^n} \big( 1 + |2^j(s-u)| \big)^{\varepsilon-2\sigma} |f_w(s)|^2  \, ds \Big\|_\infty^\frac12 \\ \!\! & \times & \!\! \Big\| \int_{|s| > |u|} (1 + |2^j(s-u)|)^{2\sigma-\varepsilon} |K_{cj}(s-u)|^2 \, ds \Big\|_\infty^\frac12 \\ \!\! & \le & \!\! C_\varepsilon \underbrace{2^{-\frac{jn}{2}} \Big\| \int_{|s| > |u|} \big(1 + |2^j s| \big)^{2\sigma-\varepsilon} |K_{cj}(s)|^2 \, ds \Big\|_\infty^\frac12}_{\mathrm{B}_{11}'(j)} \|f\|_\infty.
\end{eqnarray*}
Extracting the $\varepsilon$-power in the integral and rearranging, we get 
$$\mathrm{B}_{11}'(j) \le \big(1 + 2^j |u| \big)^{-\frac{\varepsilon}{2}} \Big\| \int_{\R^n} \big| \big(1 + |s| \big)^{\sigma} \underbrace{ 2^{-jn} K_{cj}(2^{-j} s) }_{\Phi_j(s)} \big|^2 \, ds \Big\|_\infty^\frac12.$$
This proves that
$$\mathrm{B}_{11} \, \le \ C_\varepsilon \Big( \sum_{2^j |u| > 1} (2^j |u|)^{-\frac{\varepsilon}{2}} \Big) \, \sup_{(j,x) \in \Z \times \R^n} \big\| \psi M_c(x,2^j \cdot) \big\|_{W_{2\sigma}} \, \|f\|_\infty$$
since the operator-valued function $\Phi_j$ equals 
$$\Phi_j(s) \! = \! \mathrm{diag} \Big( 2^{-jn} \big[ \psi_j M_c(x,\cdot) \big]^\vee(2^{-j} s) \Big) \! = \! \mathrm{diag} \Big( \big[ \psi M_c(x, 2^j (\cdot)) \big]^\vee(s) \Big).$$ This completes the estimate for $\mathrm{B}_{11}$ and the argument for $\mathrm{B}_{12}$ is verbatim the same.  
For $\mathrm{B}_2$, we set $A_{wj}(Q) = 2^j \big( \R^n \setminus (w - 5Q) \big)$ and rewrite it in terms of $\Phi_j$ as follows
\begin{eqnarray*}
\mathrm{B}_2(j) \!\!\! & = & \!\!\! \Big\| \int_{A_{wj}(Q)} \big( \Phi_j (s-2^ju) - \Phi_j(s) \big) f_w(2^{-j} s)  \, ds \Big\|_\infty \\ \!\!\! & = & \!\!\! \Big\| \int_{A_{wj}(Q)} \Big[ \int_0^1 \frac{d}{dr} \Phi_j(s-r 2^ju) \, dr \Big] f_w(2^{-j} s)  \, ds \Big\|_\infty \\ \!\!\! & \lesssim & \!\!\! \, \sup_{0 \le r \le 1} \, \Big\| \Big( \int_{\R^n} \big( 1 + |s-r 2^j u| \big)^{- 2\sigma} |f_w(2^{-j} s)|^2 \, ds \Big)^\frac12 \Big\|_\infty \\ 
\!\!\! & \times & \!\!\! |2^ju| \ \Big\| \Big( \int_{\R^n} \big( 1 + |s-r 2^j u| \big)^{2\sigma} \big| \nabla \Phi_j(s-r 2^j u) \big|^2 \, ds \Big)^\frac12 \Big\|_\infty \\ \!\!\! & \le & \!\!\! C_\varepsilon \, |2^j u| \, \Big\| \int_{\R^n} \big| \big( 1 + |s| \big)^{\sigma} \nabla \Phi_j(s) \big|^2 \, ds \Big\|_\infty^\frac12 \, \|f\|_\infty.
\end{eqnarray*}
To estimate the above integral, we note that \vskip-12pt
$$\widehat{\nabla \Phi}_j(\xi) \! = \! -2 \pi i \sum_{k=1}^n \mathrm{diag} \Big( \xi_k \psi(\xi) M_{c}(x,2^j \xi) \Big) \otimes e_k.$$ In particular
\begin{eqnarray*}
\lefteqn{\hskip-20pt \Big\| \int_{\R^n} \big| \big( 1 + |s| \big)^{\sigma} \nabla \Phi_j(s) \big|^2 \, ds \Big\|_\infty^\frac12} \\ \hskip20pt \!\!\! & = & \!\!\! 2 \pi \Big\| \sum_{k=1}^n \mathrm{diag} \Big( \int_{\R^n} \big| \big( 1 + |s| \big)^{\sigma} \big( \langle \hskip2pt \cdot, e_k \rangle \hskip1pt \psi M_c(x, 2^j \cdot) \big)^\wedge \big|^2 \, ds \Big) \Big\|_\infty^\frac12 \\ \!\!\! & \le & \!\!\! 2 \pi \Big( \sum_{k=1}^n \sup_{x \in \R^n} \big\| \langle \hskip2pt \cdot, e_k \rangle \hskip1pt \psi M_c(x, 2^j \cdot) \big\|_{W_{2\sigma}}^2 \Big)^\frac12 \lesssim \sup_{x \in \R^n} \big\| \psi M_c(x, 2^j \cdot) \big\|_{W_{2\sigma}}.
\end{eqnarray*}
Therefore, we obtain the following upper bound 
$$\mathrm{B}_2 \, \le \ C_\varepsilon \Big( \sum_{2^j |u| \le 1} |2^ju| \Big) \, \sup_{(j,x) \in \Z \times \R^n} \big\| \psi M_c(x,2^j \cdot) \big\|_{W_{2\sigma}} \, \|f\|_\infty.$$ This proves \eqref{Eq-HM4} in the category of Banach spaces. The cb-norm bound holds as 
well with the same proof after matrix amplification and the row term is similar. \fin  
 
\section{\large \bf Applications and discussion} \label{AppsComts}

\subsection{HMS $p$-conditions}

Given $1 < p < \infty$, we may refine Theorem $\mathrm{A}'$ with lower Sobolev regularity orders for $S_p$-boundedness as a function of $|\frac1p - \frac12|$. This natural phenomenon was originally established by Calder\'on and Torchinsky for Euclidean Fourier multipliers \cite{CT}. The argument uses complex interpolation for analytic families of operators applying Plancherel and H\"ormander-Mikhlin theorems to control the endpoints in $L_2$ and $\mathrm{H}_1$ or BMO respectively. The modifications in our case are quite simple and left to the reader. Indeed, we just need to follow their reasoning for the two terms in $\big\bracevert \hskip-2pt M \hskip-2pt \big\bracevert_{\hskip-3pt \mathrm{HMS}}$ and use Remark \ref{Rem-RowColInterpolation} for interpolation.

\begin{corollary} \label{Cor-CT} 
If $\big| \frac1p - \frac12 \big| < \frac{\delta}{n}$ and $\frac{n}{q} < \delta < \frac{n}{2}$ with $q>2$ $$\big\| S_M: S_p(\R^n) \to S_p(\R^n) \big\|_{\mathrm{cb}} \, \le \, C_p \big\bracevert \hskip-2pt M \hskip-2pt \big\bracevert_{\hskip-3pt \mathrm{HMS}_{q \delta}}.$$ 
\end{corollary}

\begin{remark}
\emph{By \cite{Graetal}, Corollary \ref{Cor-CT} is optimal since it fails when $|\frac1p - \frac12| > \frac{\delta}{n}$.}
\end{remark}

\subsection{Operator-Lipschitz functions} Set $$M_f(x,y) = \frac{f(x) - f(y)}{x-y} \quad \mbox{for} \quad \mbox{$f: \R \to \C$} \quad \mbox{Lipschitz} \quad \mbox{and} \quad \mbox{$x \neq y$}.$$ These multipliers were considered by Arazy \cite{Ar}, who conjectured $S_p$-boundedness for $1 < p < \infty$. It implies Krein's longstanding conjecture for operator-Lipschitz functions |$\|f(A) - f(B)\|_{S_p} \le C_p \|A-B\|_{S_p}$ for any pair of self-adjoint operators $A,B$ whose difference lies in the Schatten $p$-class| and it was finally confirmed in 2011 by Potapov and Sukochev \cite{PS}. It trivially follows from our main result. 

\begin{corollary}[Arazy's conjecture] \label{Cor-Arazy} 
If $1 < p < \infty$ $$\big\| S_{M_f}: S_p(\R) \to S_p(\R) \big\|_{\mathrm{cb}} \, \le \, C_p \|f\|_{\mathrm{Lip}}.$$
\end{corollary}

\dem Since $n=1$, if $f \in C^1(\R)$ Theorem A and Remark \ref{Rmk-ThmA-cb} give
$$\big\| S_{M_f} \big\|_{\mathrm{cb}(S_p(\R))} \, \le \, C \frac{p^2}{p-1} \Big( \|M_f\|_\infty + \Big\| |x-y| \Big\{ |\partial_xM_f(x,y)| + |\partial_yM_f(x,y)| \Big\} \Big\|_\infty \Big),$$ which is trivially bounded above by the Lip-norm of $f$. In general we only know that $f$ is a.e. differentiable, but we may apply Theorem A' for $M_f$ with (say) $\sigma = \frac{3}{4}$. Then, according to \eqref{Eq-HMSComparison} we see that the above inequality still holds, despite $\partial_x M_f$ and $\partial_y M_f$ are undefined in a zero-measure set. \fin 

\noindent Other forms of Arazy's conjecture ($S_p \to S_q$ with $p \neq q$) were also solved in \cite{PST}.

\subsection{Multipliers of $\alpha$-divided differences}

If $0 < \alpha < 1$, set $$M_{\alpha f}(x,y) = \frac{f(x) - f(y)}{|x-y|^\alpha} \quad \mbox{for} \quad \mbox{$f: \R \to \C$} \quad \mbox{$\alpha$-H\"older} \quad \mbox{and} \quad \mbox{$x\neq y$}.$$ The mapping properties of $S_{M_{\alpha f}}$ have remained open so far. Since Theorem $\mathrm{A}'$ just imposes $\frac12 + \varepsilon$ derivatives in one dimension, it looks there is room for improvement in Arazy's conjecture and we confirm it below. Let $\displaystyle \|f\|_{\Lambda_\alpha} = \sup_{x \neq y} |f(x) - f(y)|/|x-y|^\alpha$.

\begin{corollary}[$\alpha$-Arazy] \label{Cor-Holder} 
If $|\frac1p - \frac12| < \min \{\alpha, \frac12\}$, then 
$$\big\| S_{M_{\alpha f}}: S_p(\R) \to S_p(\R) \big\|_{\mathrm{cb}} \, \le \, C_p \|f\|_{\Lambda_\alpha}.$$ In particular, we get complete $S_p$-boundedness for $1 < p < \infty$ as long as $\alpha \ge \frac12$. 
\end{corollary}

\dem Let $\Psi_{M,j}$ be as in \eqref{Eq-SobolevNorm}. If $|s|$ is small, we easily get 
\begin{eqnarray} \label{Eq-MC1}
\omega_{1,qjy}(s) \!\!\!\! & := & \!\!\!\! \Big( \int_\R \big| \Psi_{M,j}(x+s,y) - \Psi_{M,j}(x,y) \big|^q dx \Big)^{\frac1q} \, \lesssim \, \|f\|_{\Lambda_\alpha} |s|^\alpha, \\ \label{Eq-MC2}
\omega_{2,qjx}(s) \!\!\!\! & := & \!\!\!\! \Big( \int_\R \big| \underbrace{\Psi_{M,j}(x,y+s) - \Psi_{M,j}(x,y)}_{\Theta_{M,j,s}(x,y)} \big|^q dy \Big)^{\frac1q} \, \lesssim \, \|f\|_{\Lambda_\alpha} |s|^\alpha. 
\end{eqnarray}
Indeed, if $|s|$ is small enough, the above integrals are supported by the corona-type set $\frac14 \le \min \big\{ |x-y|, |x+s-y| \big\} \le \max \big\{ |x-y|, |x+s-y| \big\} \le 4$ and we may bound them by the same integrals with $q=1$. Next, we estimate 
\begin{eqnarray*}
\big| \Theta_{M,j,s}(x,y) \big| \!\!\!\! & \le & \!\!\!\! \big| \psi(x-y-s) - \psi(x-y) \big| \, \Big| \frac{f(2^jx) - f(2^j (y+s))}{2^{j\alpha}|x-y-s|^\alpha} \Big| \\ \!\!\!\! & + & \!\!\!\! \psi(x-y) \Big| \frac{f(2^jx) - f(2^j (y+s))}{2^{j\alpha}|x-y-s|^\alpha} - \frac{f(2^j x) - f(2^j y)}{2^{j\alpha}|x-y|^\alpha} \Big| \lesssim \|f\|_{\Lambda_\alpha} |s|^\alpha.
\end{eqnarray*}
Then, if $q \ge 2$ and $0 < \delta < \alpha$, we deduce from \cite[Section V.3.5]{Stein} that
$$\big\| \Psi_{M,j}(\cdot,y) \big\|_{W_{q\delta}} \, \lesssim \, C_\varepsilon \big\| \Psi_{M,j}(\cdot,y) \big\|_q + \Big( \int_{|s| < \varepsilon} \frac{\omega_{1,qjy}(s)^2}{|s|^{1+2\delta}} ds \Big)^{\frac12} \, \le \, C_{\varepsilon \delta} \|f\|_{\Lambda_\alpha}.$$  The same bound applies to $\Psi_{M,j}(x,\cdot)$. The statement follows from Theorem $\mathrm{A'}$ with $(q,\delta)=(2,\frac{1+2\alpha}{4})$ for $\alpha > \frac12$ and from Corollary \ref{Cor-CT} for smaller values of $\alpha$. \fin

\subsection{Schur multipliers over lattices}

Now we study the form of Theorem A for Schur multipliers over integer lattices $\Z^n$. Other groups are discussed below and deeply investigated in \cite{CGPT}. Given $\phi: \Z^n \to \C$ and $j = (j_1, j_2, \ldots, j_n) \in \Z^n$ consider the discrete derivatives $$\Delta^s\phi(j) = \phi(j + e_s) - \phi(j) \quad \mbox{and}$$ \vskip-15pt 
$$\Delta^\gamma = \prod_{s=1}^n ( \Delta^s )^{\gamma_s} \quad \mbox{for higher multi-indices} \quad \gamma = (\gamma_1, \ldots, \gamma_n).$$ The operators $\Delta^\gamma$ are defined in the same way for functions on $\Z^n$ or $\R^n$. When acting on a two-variable function, we shall use $\Delta_1^\gamma = \Delta^\gamma \otimes \mathrm{id}$ and $\Delta_2^\gamma = \mathrm{id} \otimes \Delta^\gamma$ to clarify which variable does the operator act on.

\begin{corollary}[A discrete formulation] \label{Cor-Discrete} 
If $1 < p < \infty$
$$\big\| S_M \hskip-2pt: S_p(\Z^n) \to S_p(\Z^n) \big\|_{\mathrm{cb}} \, \lesssim \, \frac{p^2}{p-1} \big|\big|\big| M \big| \big| \big|_{\mathrm{HMS}_\Delta}$$ where we set $||| M |||_{\mathrm{HMS}_\Delta} := \displaystyle \sum_{|\gamma| \le [\frac{n}{2}] +1} \Big\| |j-k|^{|\gamma|} \Big\{ \big| \Delta_1^\gamma M(j,k) \big| + \big| \Delta_2^\gamma M(j,k) \big| \Big\} \Big\|_\infty.$
\end{corollary}

\dem By elementary restriction properties of Schur multipliers, it suffices from Theorem A to find a lifting function $\mathbf{M}: \R^n \times \R^n \to \C$ whose restriction to $\Z^n \times \Z^n$ coincides with $M$ and such that $$||| \mathbf{M} |||_{\mathrm{HMS}} \, \lesssim \, ||| M |||_{\mathrm{HMS}_\Delta}.$$ Lemma 4.5.1 in \cite{RT} provides $\Theta, \Phi_\gamma \in \S(\R^n)$ such that 
\begin{equation} \label{Eq-Lifting}
\sum_{k \in \Z^n} \Theta(\cdot + k) \equiv 1, \qquad \widehat{\Theta}_{\mid_{\Z^n}} \equiv \delta_0, \qquad \partial^\gamma(\widehat{\Theta}) = \overline{\Delta}^\gamma (\Phi_\gamma)
\end{equation}
for $|\gamma| \le [\frac{n}{2}]+1$ and the backwards derivative given by $$\overline{\Delta}^\gamma = \prod_{s=1}^n ( \overline{\Delta}^s )^{\gamma_s} \quad \mbox{where} \quad \overline{\Delta}^s f(\xi) = f(\xi) - f(\xi - e_s).$$ In fact, we shall use \eqref{Eq-Lifting} with $\partial^\gamma(\widehat{\Theta}) = \Delta^\gamma (\Phi_\gamma)$ instead, since the proof in \cite{RT} can be easily modified to that end. Construct the lifting multipliers by means of $\Theta$ as
$$\mathbf{M}(x,y) = \sum_{j,k \in \Z^n} \widehat{\Theta}(x-j) \widehat{\Theta}(y-k) M(j,k).$$
By the second identity in \eqref{Eq-Lifting} we clearly have $$\textbf{M} = M \quad \mbox{over} \quad \Z^n \times \Z^n.$$ Next, letting $\tau_x(f) := f(x - \cdot)$ we use the identities
$$\tau_x \circ \Delta^\gamma = (-1)^{|\gamma|} \overline{\Delta}^\gamma \circ \tau_x,$$
$$\partial^\gamma(\widehat{\Theta})(x-j) = \Delta^\gamma (\Phi_\gamma)(x-j) = (-1)^{|\gamma|} \overline{\Delta}^\gamma \big( \tau_x(\Phi_\gamma) \big) (j).$$ Then, we deduce the following inequality using summation by parts
\begin{eqnarray*}
\lefteqn{\hskip-10pt |x-y|^{|\gamma|} \big| \partial_x^\gamma \mathbf{M}(x,y) \big|} \\ [8pt] \!\! & = & \!\! |x-y|^{|\gamma|} \Big| \sum_{j,k \in \Z^n} \overline{\Delta}^\gamma (\tau_x(\Phi_\gamma))(j) \widehat{\Theta}(y-k)  M(j,k) \Big| \\ \!\! & = & \!\! |x-y|^{|\gamma|} \Big| \sum_{j,k \in \Z^n} \Phi_\gamma(x-j) \widehat{\Theta}(y-k) \Delta_1^\gamma M(j,k) \Big| \\ \!\! & \lesssim & \!\! \Big( \sum_{j,k \in \Z^n} \big| \Phi_\gamma(x-j) \widehat{\Theta}(y-k) \big| |j-k|^{|\gamma|} \big| \Delta_1^\gamma M(j,k) \big| \Big) \\
\!\! & + & \!\! \Big( \sum_{j,k \in \Z^n} |(x-j) - (y-k)|^{|\gamma|} \big| \Phi_\gamma(x-j) \widehat{\Theta}(y-k) \big| \big| \Delta_1^\gamma M(j,k) \big| \Big) \\ \!\! & \le & \!\! ||| M |||_{\mathrm{HMS}_\Delta} \Big( \sum_{j,k \in \Z^n} \big[ 1 + |(x-j) - (y-k)|^{|\gamma|} \big] \big| \Phi_\gamma(x-j) \widehat{\Theta}(y-k) \big| \Big). 
\end{eqnarray*}
The last sum above in uniformly bounded in $x,y$ since $\Phi_\gamma, \Theta$ are Schwartz functions. The derivatives of $\mathbf{M}$ in the $y$-variable are dealt with in the same way. \fin

\begin{remark}
\emph{After this paper was finished, Tao Mei pointed us he already knew Corollary \ref{Cor-Discrete} for $n=1$ controlling the dyadic 1-variation of $M_r$ and $M_c$ respectively}.
\end{remark}

\subsection{A matrix Littlewood-Paley theorem}

Let $\Psi: \R^n \times \R^n \to \R_+$ be a smooth function which generates a Littlewood-Paley partition of unity. By that we just mean that the dilated functions $\Psi_j(x,y) = \Psi(2^j x, 2^j y)$ satisfy the following conditions:
\begin{itemize}
\item[a)] $\sum_{j \in \Z} \Psi_j^2 = 1$ a.e.

\item[b)] The supports of $\Psi_j$ have finite overlapping.
\end{itemize}

\begin{corollary} [Matrix LP theorem] \label{Cor-LP} Let $1 < p < \infty:$ 
\begin{itemize}
\item[i)] If $p \le 2$ $$\|A\|_{S_p} \asymp_{\mathrm{cb}} \inf_{S_{\Psi_j}(A) = A_j+B_j} \Big\| \Big( \sum_{j \in \Z} A_jA_j^* + B_j^*B_j \Big)^{\frac12} \Big\|_{S_p}.$$

\item[ii)] If $p \ge 2$ $$\|A\|_{S_p} \asymp_{\mathrm{cb}} \Big\| \Big( \sum_{j \in \Z} S_{\Psi_j}(A)S_{\Psi_j}(A)^* + S_{\Psi_j}(A)^*S_{\Psi_j}(A) \Big)^{\frac12} \Big\|_{S_p}.$$
\end{itemize}
The constants in the cb-isomorphisms are comparable to $p^2/(p-1)$ as $p \to 1,\infty$.
\end{corollary}

\dem The terms on the right hand side are usually denoted in the literature as the $S_p[RC_p]$-norm of the sequence $(S_{\Psi_j}(A))_{j \in \Z}$. The norm takes values in the operator spaces $R_p+C_p$ or $R_p \cap C_p$, according to the value of $p$ as above. By the noncommutative Khintchine inequality \cite{LP}, both mostright terms may be written as the $\pm 1$ conditional expectation $$\mathbf{E}_\varepsilon \Big\| \sum_{j \in \Z} \varepsilon_j S_{\Psi_j}(A) \Big\|_{S_p} \le_{\mathrm{cb}} C_p \, \mathbf{E}_\varepsilon \Big|\Big|\Big| \sum_{j \in \Z} \varepsilon_j \Psi_j \Big|\Big|\Big|_{\mathrm{HMS}} \|A\|_{S_p}.$$ The above inequality follows from Theorem A. Property b) of $\Psi_j$ gives
$$\mathbf{E}_\varepsilon \Big|\Big|\Big| \sum_{j \in \Z} \varepsilon_j \Psi_j \Big|\Big|\Big|_{\mathrm{HMS}} \, \lesssim \, \sup_{j \in \Z} |||\Psi_j|||_{\mathrm{HMS}} \, = \, ||| \Psi |||_{\mathrm{HMS}} \, < \, \infty.$$
This proves the lower bounds. The upper bounds follow from a) and duality
$$\|A\|_{S_p} \, = \sup_{\|B\|_{S_{p'}} \le 1} \sum_{j \in \Z} \mathrm{tr} \big( S_{\Psi_j}(A) S_{\Psi_j}(B)^* \big) \, \lesssim_{\mathrm{cb}} \, \big\|(S_{\Psi_j}(A)) \big\|_{S_p[RC_p]},$$ using the well-known fact that $S_p[RC_p]^* = S_{p'}[RC_{p'}]$. Recall that we only use one side of the Khintchine inequalities which holds up to constants independent of $p$. Therefore, the final dependence on $p$ is only dictated by Theorem A. \fin

\begin{remark}
\emph{Let $\psi_0: \R_+ \to \R_+$ denote the radial generator of the function $\psi$ constructed in \eqref{Eq-LP}. The choice $\Psi(x,y) = \psi_0(|x-y|)$ gives the Toeplitz form of Corollary \ref{Cor-LP}, which by Fourier-Schur transference is known to be equivalent to the classical Littlewood-Paley theorem. What appears to be new is the following choice $$\Psi(x,y) = \psi_0 \Big\{ \big( |x|^2 + |y|^2 \big)^\frac12 \Big\},$$ giving smooth coronas $S_{\Psi_j}(A) = \mathrm{Cor}_j(A)$ with nonflat boundary. If $p \ge 2$ we get $$\|A\|_{S_p} \asymp_{\mathrm{cb}} \Big\| \Big( \sum_{j \in \Z} \mathrm{Cor}_j(A) \mathrm{Cor}_j(A)^* + \mathrm{Cor}_j(A)^* \mathrm{Cor}_j(A) \Big)^{\frac12} \Big\|_{S_p}.$$}
\end{remark}

\subsection{Fourier multipliers in group algebras}

Let $\beta: \G \to \R^n$ be an orthogonal cocycle and assume that $m: \G \to \C$ satisfies $m = \widetilde{m} \circ \beta$. The main discovery in \cite{GJP,JMP1,JMP2} was that a HM theory in group von Neumann algebras is possible in terms of the $\beta$-lifted symbols $\widetilde{m}$. The main drawbacks of this theory are that it requires orthogonal actions, regularity is measured in terms of the cocycle dimension (which excludes infinite-dimensional cocycles) and imposes auxiliary differential structures which appear to be unnecessary or at least less natural for Lie groups. Under these difficulties, the works \cite{MR,MRX,PRS} provide partial answers for the free group and the special linear groups $S \hskip-1pt L_n(\R)$. However, very little is known on optimal regularity conditions for $L_p$-boundedness of Fourier multipliers on Lie group algebras $\V$ in terms of canonical metrics and Lie differentiation. More precisely, let us write $d_g^\gamma m(g)$ for the left-invariant Lie derivative of order $\gamma = (j_1, j_2, \ldots, j_{|\gamma|})$ $(1 \le j_i \le n)$ of  the Fourier symbol $m$ at $g$. If $\G$ is a unimodular Lie group with $\dim \G = n$ we want to investigate natural metrics $L: \G \times \G \to \R_+$ for which the following inequality holds for $1 < p < \infty$ 
\begin{equation} \tag{HM$_\G$} \label{Eq-HM}
\big\| T_m: L_p(\V) \to L_p(\V) \big\|_{\mathrm{cb}} \, \lesssim \, C_p \sum_{|\gamma| \le [\frac{n}{2}]+1} \big\| L(g,e)^{|\gamma|} d_g^\gamma m(g) \big\|_\infty.
\end{equation}
Our first result is a \textbf{local HM theorem}: Inequality \eqref{Eq-HM} holds for arbitrary Lie groups with the Riemannian metric for symbols with small support. According to \cite{LdlS}, locality is necessary for such a general statement. Next, a \textbf{HM theorem for stratified Lie groups}:  Inequality \eqref{Eq-HM} holds for the subRiemannian metric and arbitrary symbols. This uses a stratifed Mikhlin condition (a derivative in the $k$-th stratum must be dealt with as a $k$-th order derivative) which is substantially and necessarily different from the (dual) spectral approach developed by Cowling, M\"uller, Ricci, Stein, etc... Finally, a \textbf{HM theorem for high rank simple Lie groups}: Inequality \eqref{Eq-HM} holds for general symbols and a locally Euclidean metric whose assymptotic behavior is dictated by the adjoint representation. This goes far beyond the nearly optimal results for $\SL$ in \cite{PRS}. 

These results have been the object of our desire for quite some time. Among other ingredients, the proof requires a generalized form of Theorem A which we postpone to \cite{CGPT}. The Introduction there provides historical precedents, a detailed account of the motivations and rigorous statements of these $L_p$-inequalities. Such form of Theorem A also expands the main multiplier theorems in \cite{JMP1,JMP2}.

\noindent \textbf{Acknowledgement.} JP wants to express his gratitude for very valuable comments on preliminary versions to Tao Mei, \'Eric Ricard in several occasions and Fedor Sukochev. We are also indebted to the referees for a very careful reading. All the authors were supported by the Spanish Grant PID2019-107914GB-I00 \lq\lq Fronteras del An\'alisis Arm\'onico" (MCIN / PI J. Parcet) as well as Severo Ochoa Grant CEX2019-000904-S (ICMAT), funded by MCIN/AEI 10.13039/501100011033. Jos\'e Conde-Alonso was also supported by the Madrid Government Program V PRICIT Grant SI1/PJI/2019-00514. Eduardo Tablate was supported as well by Spanish Ministry of Universities with a FPU Grant with reference FPU19/00837.    

\bibliographystyle{amsplain}

%\newpage

\enlargethispage{2,5cm}

\vskip4pt

\small

\noindent \textbf{Jos\'e M. Conde-Alonso} \hfill \textbf{Adri\'an M. Gonz\'alez-P\'erez} \\ 
Universidad Aut\'onoma de Madrid \hfill Universidad Aut\'onoma de Madrid
\\ Instituto de Ciencias Matem{\'a}ticas \hfill Instituto de Ciencias Matem{\'a}ticas
%\\ \null \hfill C/ Nicol\'as Cabrera 13-15. 28049, Madrid. Spain 
\\ \texttt{jose.conde@uam.es} \hfill \texttt{adrian.gonzalez@uam.es}

\vskip0pt

\noindent \textbf{Javier Parcet} \hfill \noindent \textbf{Eduardo Tablate} \\
Instituto de Ciencias Matem{\'a}ticas \hfill Instituto de Ciencias Matem{\'a}ticas 
\\ CSIC \hfill CSIC %Universidad Aut\'onoma de Madrid  
%\\ \null \hfill C/ Nicol\'as Cabrera 13-15. 28049, Madrid. Spain 
\\ \texttt{parcet@icmat.es} \hfill \texttt{eduardo.tablate@icmat.es}

%\hfill \noindent \textbf{Javier Parcet} \\
%\null \hfill Instituto de Ciencias Matem{\'a}ticas 
%\\ \null \hfill Consejo Superior de Investigaciones Cient{\'\i}ficas 
%\\ \null \hfill C/ Nicol\'as Cabrera 13-15. 28049, Madrid. Spain 
%\\ \null \hfill\texttt{parcet@icmat.es}

%\vskip2pt

%\hfill \noindent \textbf{Eduardo Tablate} \\
%\null \hfill Instituto de Ciencias Matem{\'a}ticas 
%\\ \null \hfill Universidad Aut\'onoma de Madrid 
%\\ \null \hfill C/ Nicol\'as Cabrera 13-15. 28049, Madrid. Spain 
%\\ \null \hfill\texttt{eduardo.tablate@estudiante.uam.es}

\end{document}